\documentclass[graybox, envcountsame,envcountsect]{svmult}

\usepackage{helvet}          
\usepackage{courier}         
\usepackage{amsmath}
\usepackage{amssymb}
\usepackage[all]{xy}
\usepackage{comment} 
\usepackage{amscd,verbatim, color}
\usepackage{enumerate}  
\usepackage{amsbsy} 
\usepackage[mathscr]{eucal}
\usepackage{mathrsfs} 
\usepackage{array}
\usepackage{rotating} 
\usepackage{chngcntr} 

\smartqed
\newenvironment{qedproof}{\begin{proof}}{\qed\end{proof}}  

\spnewtheorem{rem}[theorem]{Remark}{\bfseries}{}
\spnewtheorem{examp}[theorem]{Example}{\bfseries}{}
\spnewtheorem{conj}[theorem]{Conjecture}{\bfseries}{\itshape}
\spnewtheorem{exmp}[theorem]{Example}{\bfseries}{}
\spnewtheorem{warning}[theorem]{Warning}{\bfseries}{}
\spnewtheorem{ques}[theorem]{Question}{\bfseries}{\normalfont}

\spnewtheorem{fact}[theorem]{Fact}{\bfseries}{\itshape}
\spnewtheorem{observation}[theorem]{Observation}{\bfseries}{\itshape}
\spnewtheorem{false}[theorem]{False Statement}{\bfseries}{\itshape}
\spnewtheorem{counterexample}[theorem]{Counterexample}{\bfseries}{\normalfont}
\spnewtheorem{def-lemma}[theorem]{Definition-Lemma}{\bfseries}{\itshape}

\spnewtheorem*{notation}{Notation}{\bfseries}{\normalfont}  
\spnewtheorem*{algorithm}{Algorithm}{\bfseries}{\rmfamily}

\newcommand{\dual}[2]{\widehat{#1}_{\operatorname{#2}}}
\newcommand{\invHom}[3]{\operatorname{Hom}_{#1}({#2},{#3})}

\begin{document}
\motto{\flushright \normalsize Dedicated to David Vogan on the occasion of his 60th birthday, with admiration of his epoch-making contributions to the field}
\title*{A program for branching problems in the representation theory of real reductive groups}
\author{Toshiyuki Kobayashi}
\institute{Toshiyuki Kobayashi 
\at Kavli IPMU (WPI)
 and Graduate School of Mathematical Sciences, 
The University of Tokyo
\\ \email{toshi@ms.u-tokyo.ac.jp}}
\maketitle

\counterwithin{equation}{section}   
\renewcommand{\theequation}{\thesection.\arabic{equation}}  

\abstract{We wish to understand 
how irreducible representations
 of a group $G$ behave 
 when restricted to a subgroup $G'$
(the {\it{branching problem}}).  
Our primary concern is with representations
 of reductive Lie groups,
 which involve both algebraic and analytic approaches.  
We divide branching problems
 into three stages:
 (A) abstract features
 of the restriction;
 (B) branching laws (irreducible decompositions
     of the restriction);
 and (C) construction of symmetry breaking operators
 on geometric models.  
We could expect a simple and detailed study
 of branching problems in Stages B and C
 in the settings
 that are {\it{a priori}} known
 to be \lq\lq{nice}\rq\rq\
 in Stage A, 
and conversely,
 new results and methods
 in Stage C that
 might open another fruitful direction
 of branching problems
 including Stage A.  
The aim of this article
 is to give new perspectives
 on the subjects, 
 to explain the methods 
 based on some recent progress, 
 and to raise some conjectures
 and open questions.}
\keywords
{branching law, 
symmetry breaking operator, 
unitary representation, 
Zuckerman--Vogan's $A_{{\mathfrak{q}}}(\lambda)$ module,
reductive group,
 spherical variety, 
multiplicity-free representation \\ \\
\textbf{MSC (2010):} Primary 22E46; Secondary 53C35}%

\section{Program --- ABC for branching problems}

{}From the viewpoint
 of analysis and synthesis, 
one of the fundamental problems
 in representation theory 
 is to classify the smallest objects
 (e.g., 
 irreducible representations), 
and another is to understand
 how a given representation can be built up from the smallest objects
 (e.g., 
 irreducible decomposition).  
A typical example of the latter
 is the {\it{branching problem}}, 
 by which we mean the problem of understanding 
 how irreducible representations $\pi$
 of a group $G$ behave 
 when restricted to subgroups $G'$.  
We write $\pi|_{G'}$ 
 for a representation $\pi$
 regarded as a representation of $G'$.  
Our primary concern is with real reductive Lie groups.  
We propose a program
 for branching problems
 in the following three stages:
\vskip 0.8pc
\par\noindent
{\bf{Stage A.}}\enspace
Abstract features of the restriction $\pi|_{G'}$.  
\par\noindent
{\bf{Stage B.}}\enspace
Branching laws.  
\par\noindent
{\bf{Stage C.}}\enspace
Construction of symmetry breaking operators.  
\vskip 0.8pc
Here, 
 by a {\it{symmetry breaking operator}}
 we mean a continuous $G'$-homomorphism from 
 the representation space of $\pi$
 to that of an irreducible representation $\tau$ of the subgroup $G'$.

\par
Branching problems
 for infinite-dimensional representations
 of real reductive groups 
 involve various aspects.  
Stage A involves several aspects
 of the branching problem, 
 among which we highlight that of multiplicity
 and spectrum here:

\vskip 0.5pc
{\bf{A.1.}}\enspace
Estimates of multiplicities of irreducible representations
 of $G'$ occurring in the restriction $\pi|_{G'}$
of an irreducible representation
 $\pi$ of $G$.  
 (There are several \lq\lq{natural}\rq\rq\ but inequivalent definitions
 of multiplicities, 
 see Sections \ref{sec:irreddeco} and \ref{subsec:2}.)
Note that:
\begin{enumerate}
\item[$\bullet$]
multiplicities of the restriction $\pi|_{G'}$
may be infinite even when $G'$ is a maximal subgroup in $G$;
\item[$\bullet$]
multiplicities may be at most one
 (e.g., 
 Howe's theta correspondence \cite{xhowe}, 
 Gross--Prasad conjecture \cite{GP}, 
 visible actions \cite{mfbundl}, 
 etc.).  
\end{enumerate}

\vskip 0.5pc
{\bf{A.2.}}\enspace
Spectrum of the restriction $\pi|_{G'}$:
\begin{enumerate}
\item[$\bullet$]
(discretely decomposable case)\enspace
branching problems may be
 purely algebraic and combinatorial
 (\cite{xdv, E-W, xgrwaII, xkInvent94, xkAnn98, xkInvent98, xkdecoaspm, decoAq, xkyosh13, xspeh});
\item[$\bullet$]
(continuous spectrum)\enspace
branching problems may have analytic features
 \cite{CKOP, xtsbon, Mol, Re79}.  
(For example, 
 some special cases
 of branching laws
 of unitary representations are equivalent
 to a Plancherel-type theorem
 for homogeneous spaces.)  
\end{enumerate}
The goal of Stage A in branching problems 
 is to analyze the aspects
 A.1 and A.2
 in complete generality.  
A theorem in Stage A would 
 be interesting
 on its own, 
but might also serve as
 a foundation for further detailed study
 of the restriction $\pi|_{G'}$
 (Stages B and C).  
An answer in Stage A may also 
 suggest an approach 
 depending on specific features
 of the restrictions.  
For instance,
 if we know {\it{a priori}}
 that the restriction 
 $\pi|_{G'}$ is discretely decomposable 
 in Stage A, 
then one might use algebraic methods
 (e.g., 
 combinatorics,
 ${\mathcal{D}}$-modules,
 etc.)
 to attack Stage B.  
If the restriction $\pi|_{G'}$
 is known {\it{a priori}}
 to be multiplicity-free
 in Stage A, 
 one might expect to find 
 not only explicit irreducible decompositions
 (Stage B)
 but also quantitative estimates 
 such as $L^p-L^q$ estimates, 
 and Parseval--Plancherel type theorems 
 for branching laws 
 (Stage C).

In this article, 
 we give some perspectives
of the subject based on a general theory on A.1 and A.2, 
 and recent progress in some classification theory:
\begin{enumerate}
\item[$\bullet$]
the multiplicities to be finite
 [bounded, one, $\cdots$], 
\item[$\bullet$]
the spectrum to be discrete / continuous.  
\end{enumerate}
We also discuss a new phenomenon
 ({\it{localness theorem}}, Theorem \ref{thm:local})
 and open questions.

Stage B concerns the irreducible decomposition
 of the restriction.  
For a finite-dimensional representation
 such that the restriction $\pi|_{G'}$ is completely reducible, 
 there is no ambiguity
 on a meaning of the irreducible decomposition.  
For a unitary representation $\pi$,
 we can consider Stage B
 by using the direct integral
 of Hilbert spaces (Fact \ref{fact:irreddeco}).  
However,
 we would like to treat 
 a more general setting 
 where $\pi$ is not necessarily
 a unitary representation.  
In this case, 
 we may consider Stage B
 as the study of 
\begin{equation}
\label{eqn:order}
  \invHom{G'}{\pi|_{G'}}{\tau}
  \quad
  \text{or}
  \quad
  \invHom{G'}{\tau}{\pi|_{G'}}
\end{equation}
 for irreducible representations
 $\pi$ and $\tau$ of $G$ and $G'$, 
respectively.

Stage C is more involved
 than Stage B 
 as it asks for concrete intertwining operators
 (e.g., 
 the projection operator
 to an irreducible summand)
 rather than an abstract decomposition;
it asks for the decomposition
 of vectors
 in addition to that of representations.  
Since Stage C depends
 on the realizations
 of the representations; 
 it often interacts with geometric 
 and analytic problems.

We organize
 this article 
 not in the natural order, 
 Stage A $\Rightarrow$ Stage B $\Rightarrow$ Stage C, 
 but in an opposite order, 
Stage C $\Rightarrow$ Stages A and B.  
This is because it is only recently 
 that a complete construction 
of all symmetry breaking operators
 has been carried out
 in some special settings, 
 and because such examples 
 and new methods 
 might yield yet
 another interesting direction 
 of branching problems in Stages A to C.  
The two spaces in \eqref{eqn:order}
 are discussed 
in Sections \ref{sec:ep}--\ref{sec:A2} from 
different perspectives
 (Stage A).  
The last section returns to Stage C
 together with comments on the general theory 
 (Stages A and B).

\section{Two concrete examples of Stage C}
\label{sec:ex}
In this section, 
 we illustrate Stage C
 in the branching program
 with two recent examples, 
namely,
 an explicit construction
 and a complete classification
 of {\it{differential}} symmetry breaking
 operators
 (Section \ref{subsec:SL2})
 and {\it{continuous}} symmetry breaking operators
 (Section \ref{subsec:2.2}).    
They have been carried out only 
 in quite special situations
 until now.  
In this section 
 we examine these new examples
 by making some observations
 that may contain some interesting hints
 for future study.  
In later sections,
 we discuss to what extent the new results
 and methods
 apply to other situations
 and what the limitations of the general theory for Stage A
 would be.  

\subsection{Rankin--Cohen bidifferential operators for the tensor products of ${\mathrm{SL}}_2$-modules}
\label{subsec:SL2}

Taking the ${\mathrm{SL}}_2$-case
 as a prototype,
 we explain 
 what we have in mind
 for Stage C
 by comparing it with Stages A and B.  
We focus on {\it{differential}} symmetry breaking operators
 in this subsection, 
 and point out that there are some missing operators
 even in the classical ${\mathrm{SL}}_2$-case
 (\cite{C75, xrankin}, 
  see also van Dijk--Pevzner \cite{vDijkPevzner}, 
Zagier \cite{xzagier94}).

First, 
 we begin with finite-dimensional representations.  
For every $m \in {\mathbb{N}}$, 
there exists
 the unique $(m+1)$-dimensional irreducible holomorphic representation
 of ${\mathrm{SL}}(2,{\mathbb{C}})$.  
These representations can be realized on the space
 $\operatorname{Pol}_m[z]$
 of polynomials in $z$
 of degree at most $m$, 
 by the following action of ${\mathrm{SL}}(2,{\mathbb{C}})$
 with $\lambda=-m$:
\begin{equation}
\label{eqn:czd}
   (\pi_{\lambda}(g) f)(z)
   =
    (cz+d)^{-\lambda}f(\frac{az+b}{cz+d})
\quad
  \text{for }
\,\,
  g^{-1}=\begin{pmatrix} a & b \\ c & d \end{pmatrix}.  
\end{equation}
The tensor product of two such representations
 decomposes into irreducible representations of ${\mathrm{SL}}(2,{\mathbb{C}})$
 subject to the classical Clebsch--Gordan formula:

\begin{equation}
\label{eqn:CG}
   \operatorname{Pol}_{m}[z] \otimes \operatorname{Pol}_{n}[z]
    \simeq
   \operatorname{Pol}_{m+n}[z] \oplus \operatorname{Pol}_{m+n-2}[z] \oplus
    \cdots
    \oplus
    \operatorname{Pol}_{|m-n|}[z].  
\end{equation}

Secondly,
 we recall an analogous result
 for infinite-dimensional representations
 of ${\mathrm{SL}}(2,{\mathbb{R}})$.  
For this,
 let $H_+$ be the Poincar{\'e} upper half plane
 $\{z \in {\mathbb{C}}: \operatorname{Im}z >0\}$.  
Then ${\mathrm{SL}}(2,{\mathbb{R}})$ acts on the space
 ${\mathcal{O}}(H_+)$
 of holomorphic functions on $H_+$
 via $\pi_{\lambda}$ ($\lambda \in {\mathbb{Z}}$).  
Further,
 we obtain an irreducible unitary representation of ${\mathrm{SL}}(2,{\mathbb{R}})$
 on the following Hilbert space
 $V_{\lambda}$
 (the {\textit{weighted Bergman space}})
 via $\pi_{\lambda}$ for $\lambda >1$:
\[
   V_{\lambda}:=
   \{ f \in {\mathcal{O}}(H_+):
   \int_{H_+} |f(x+\sqrt{-1}y)|^2 y^{\lambda-2} dx dy < \infty
\}, 
\]
where the inner product is given by
\[
   \int_{H_+} f(x+\sqrt{-1}y) \overline{g(x+\sqrt{-1}y)} y^{\lambda-2} dx dy
\quad
\text{for $f, g \in V_{\lambda}$}.  
\]

Repka \cite{Re79} and Molchanov \cite{Mol}
 obtained the irreducible decomposition
 of the tensor product of two such unitary representations, 
 namely, 
 there is a unitary equivalence
 between unitary representations
 of ${\mathrm{SL}}(2,{\mathbb{R}})$:
\begin{equation}
\label{eqn:Repka}
   V_{\lambda_1} \widehat \otimes V_{\lambda_2}
   \simeq
{\sum_{a=0}^{\infty}}{}^{\oplus}V_{\lambda_1+ \lambda_2+ 2a}, 
\end{equation}
where the symbols $\widehat \otimes$ and $\sum^{\oplus}$ denote 
 the Hilbert completion
 of the tensor product $\otimes$ 
 and the algebraic direct sum $\oplus$, 
respectively.  
We then have:
\begin{observation}
\begin{enumerate}[{\rm (1)}] 
\item {\rm{(multiplicity)}}
\enspace
Both of the irreducible decompositions \eqref{eqn:CG}
 and \eqref{eqn:Repka}
 are multiplicity-free.  
\item 
{\rm{(spectrum)}}\enspace
There is no continuous spectrum
 in either of the decompositions 
 \eqref{eqn:CG}
 or \eqref{eqn:Repka}.  
 \end{enumerate}
\end{observation}
These abstract features (Stage A) are immediate consequences
 of the decomposition formul{\ae} \eqref{eqn:CG} and \eqref{eqn:Repka}
 (Stage B), 
however,
 one could tell 
 these properties 
without explicit formul{\ae} from the general theory
 of visible actions
 on complex manifolds
 \cite{xrims40, mfbundl}
 and a general theory
 of discrete decomposability
 \cite{xkInvent94, xkAnn98}.  
For instance, 
the following holds:

\begin{fact}
\label{fact:2.2}
Let $\pi$ be an irreducible unitary highest weight representation 
 of a real reductive Lie group $G$, 
 and $G'$ a reductive subgroup of $G$.  
\begin{enumerate}[{\rm (1)}]
\item
{\rm{(multiplicity-free decomposition)}}\enspace
The restriction $\pi|_{G'}$ is multiplicity-free
 if $\pi$ has a scalar minimal $K$-type
 and 
 $(G,G')$ is a symmetric pair.  
\item
{\rm{(spectrum)}}
The restriction $\pi|_{G'}$ is discretely decomposable
 if the associated Riemannian symmetric spaces
 $G/K$ and $G'/K'$
 carry Hermitian symmetric structures
 such that the embedding
 $G'/K' \hookrightarrow G/K$ is holomorphic.  
\end{enumerate}
\end{fact}

Stage C asks for a construction of the following explicit ${\mathrm{SL}}_2$-intertwining
 operators ({\it{symmetry breaking operators}}):
\begin{alignat*}{2}
   \operatorname{Pol}_{m}[z] \otimes \operatorname{Pol}_{n}[z]
   &\to 
   \operatorname{Pol}_{m+n-2a}[z] 
\qquad
&&\text{for $0 \le a \le \operatorname{min}(m,n)$, }
\\
   V_{\lambda_1} \widehat \otimes V_{\lambda_2}
   &\to 
   V_{\lambda_1+ \lambda_2+ 2a}
\qquad
&&\text{for $a \in {\mathbb{N}}$},   
\end{alignat*}
 for finite-dimensional and infinite-dimensional representations,
 respectively.  
We know a priori from Stages A and B
 that such intertwining operators exist uniquely 
 (up to scalar multiplications)
 by Schur's lemma 
 in this setting.  
A (partial) answer to this question 
 is given by the classical Rankin--Cohen bidifferential operator, 
 which is defined by
\begin{multline*}
{\mathcal{RC}}_{\lambda_1,\lambda_2}^{\lambda_1+\lambda_2+2a}(f_1,f_2)(z)
\\
   :=
\sum_{l=0}^a
  \frac{(-1)^l}{l ! (a-l)!}
  \frac{(\lambda_1+a-1) ! (\lambda_2+a-1)!}{(\lambda_1+a-l-1)! (\lambda_2+l-1)!}  \frac{\partial^{a-l}f_1}{\partial z^{a-l}}(z)
  \frac{\partial^{l}f_2}{\partial z^{l}}(z)
\end{multline*}
for $a \in {\mathbb{N}}$, 
 $\lambda_1, \lambda_2 \in \{2,3,4,\dots\}$, 
 and $f_1$, $f_2 \in {\mathcal{O}}(H_+)$.  
Then 
$
{\mathcal{RC}}_{\lambda_1,\lambda_2}^{\lambda_1+\lambda_2+2a}
$
 is an operator which intertwines 
 $\pi_{\lambda_1} \widehat \otimes \pi_{\lambda_2}$
 and $\pi_{\lambda_1+\lambda_2+2a}$.

More generally, 
 we treat {\it{non-unitary}}
 representations $\pi_{\lambda}$
 on ${\mathcal{O}}(H_+)$
 of the universal covering group ${\mathrm{SL}}(2,{\mathbb{R}})^{\sim}$
 of ${\mathrm{SL}}(2,{\mathbb{R}})$ 
 by the same formula \eqref{eqn:czd} for $\lambda \in {\mathbb{C}}$,
 and consider a continuous linear map
\begin{equation}
\label{eqn:SL2T}
   T: {\mathcal{O}}(H_+ \times H_+) \to {\mathcal{O}}(H_+)
\end{equation}
 that intertwines $\pi_{\lambda_1} \otimes \pi_{\lambda_2}$
 and $\pi_{\lambda_3}$, 
 where ${\mathrm{SL}}(2,{\mathbb{R}})^{\sim}$ 
 acts on ${\mathcal{O}}(H_+ \times H_+)$
 via $\pi_{\lambda_1} \otimes \pi_{\lambda_2}$
 under the diagonal action.  
We denote by $H(\lambda_1, \lambda_2, \lambda_3)$
the vector space 
 of symmetry breaking operators $T$
 as in \eqref{eqn:SL2T}.  
\begin{ques}
\label{que:SL2}
\begin{enumerate}[{\rm (1)}]
\item
{\rm{(Stage B)}}
Find the dimension of $H(\lambda_1, \lambda_2, \lambda_3)$
 for $(\lambda_1, \lambda_2, \lambda_3) \in {\mathbb{C}}^3$.  
\item
{\rm{(Stage C)}}
Explicitly construct a basis of $H(\lambda_1, \lambda_2, \lambda_3)$
 when it is nonzero.  
\end{enumerate}
\end{ques}
Even in the ${\mathrm{SL}}_2$-setting, 
we could not find a complete answer to Question \ref{que:SL2}
 in the literature, 
 and thus we explain our solution below.

Replacing $\mu !$ by $\Gamma(\mu+1)$, 
 we can define 
\begin{equation}
\label{eqn:RC}
   {\mathcal{RC}}_{\lambda_1,\lambda_2}^{\lambda_3}(f_1,f_2)(z)
   :=
  \sum_{l=0}^a
  \frac{(-1)^l}{l ! (a-l)!}
  \frac{\Gamma(\lambda_1+a)\Gamma(\lambda_2+a)}{\Gamma(\lambda_1+a-l)\Gamma(\lambda_2+l)}
  \frac{\partial^{a-l}f_1}{\partial z^{a-l}}(z)
  \frac{\partial^{l}f_2}{\partial z^{l}}(z), 
\end{equation}
where $a:=\frac 1 2 (\lambda_3 - \lambda_1 - \lambda_2)$
 as long as $(\lambda_1, \lambda_2, \lambda_3)$
 belongs to 
\[
\Omega := \{(\lambda_1, \lambda_2, \lambda_3) \in {\mathbb{C}}^3:
     \lambda_3-\lambda_1-\lambda_2 = 0,2,4, \dots\}.  
\]

We define a subset 
 $\Omega_{\operatorname{sing}}$ 
 of $\Omega$ by 
\[
\Omega_{\operatorname{sing}}
:=
\{
(\lambda_1, \lambda_2, \lambda_3) \in \Omega
:
   \lambda_1, \lambda_2, \lambda_3 \in {\mathbb{Z}}, 
 \quad
\lambda_3-|\lambda_1 - \lambda_2| \ge 2 \ge \lambda_1 + \lambda_2 
+ \lambda_3
\}. 
\]
Then we have the following classification 
 of symmetry breaking operators 
 by using the \lq\lq{F-method}\rq\rq\
 ({\cite[Part II]{KP1}}).  
Surprisingly,
 it turns out
 that any symmetry breaking operator
 \eqref{eqn:SL2T} is given by a differential operator.  

\begin{theorem}
\label{thm:RC}
\begin{enumerate}[{\rm (1)}]
\item  
$H(\lambda_1, \lambda_2, \lambda_3) \ne \{0\}$
 if and only if $(\lambda_1, \lambda_2, \lambda_3) \in \Omega$.
\vskip 0.5pc
{}  From now on, 
 we assume $(\lambda_1, \lambda_2, \lambda_3) \in \Omega$.  
\item  
$\dim_{\mathbb{C}} H(\lambda_1, \lambda_2, \lambda_3)=1$
 if and only if ${\mathcal{RC}}_{\lambda_1,\lambda_2}^{\lambda_3} \ne 0$, 
 or equivalently, 
 
 $(\lambda_1, \lambda_2, \lambda_3) \not\in \Omega_{\operatorname{sing}}$.  
In this case, 
 $H(\lambda_1, \lambda_2, \lambda_3)={\mathbb{C}} \ {\mathcal{RC}}_{\lambda_1,\lambda_2}^{\lambda_3}$.  
\item 
The following three conditions
 on $(\lambda_1, \lambda_2, \lambda_3) \in \Omega$
 are equivalent:
\begin{enumerate}[{\rm (i)}]
\item 
$\dim_{\mathbb{C}} H(\lambda_1, \lambda_2, \lambda_3)=2$.  
\item 
${\mathcal{RC}}_{\lambda_1,\lambda_2}^{\lambda_3} =0$.  
\item 
$(\lambda_1, \lambda_2, \lambda_3) \in \Omega_{\operatorname{sing}}$.  
\end{enumerate}
In this case,
 the two-dimensional vector space
 $H(\lambda_1, \lambda_2, \lambda_3)$
is spanned by 
\[
  {\mathcal{RC}}_{2-\lambda_1,\lambda_2}^{\lambda_3} 
   \circ 
   ((\frac{\partial}{\partial z_1})^{1- \lambda_1} \otimes {\operatorname{id}})\quad
\text{ and }
\quad
{\mathcal{RC}}_{\lambda_1,2-\lambda_2}^{\lambda_3}
\circ 
   (\operatorname{id} \otimes (\frac{\partial}{\partial z_2})^{1- \lambda_2}).
\]
\end{enumerate}
\end{theorem}
Theorem \ref{thm:RC} answers Question \ref{que:SL2} (1) and (2).  
Here are some observations.

\begin{observation}
\label{obs:SL2}
\begin{enumerate}[{\rm (1)}]
\item 
{\rm{(localness property)}}\enspace
Any symmetry breaking operator from 
 $\pi_{\lambda_1} \otimes \pi_{\lambda_2}$ to 
 $\pi_{\lambda_3}$ is given 
 by a differential operator
 in the holomorphic realization 
of $\pi_{\lambda_j}$ ($j=1,2,3$).  
\item 
{\rm{(multiplicity-two phenomenon)}}\enspace
The dimension of the space
 of symmetry breaking operators
 jumps up exactly 
 when the holomorphic continuation 
 of the Rankin--Cohen bidifferential operator vanishes.  
\end{enumerate}
\end{observation}

The localness property
 in Observation \ref{obs:SL2} (1)
 was recently proved 
 in a more general setting 
 (see Theorem \ref{thm:local} and Conjecture \ref{conj:local}).  

\begin{rem}
[higher multiplicities at $\Omega_{\operatorname{sing}}$]
\label{rem:2.6}
\begin{enumerate}[{\rm (1)}] 
\item 
{} From the viewpoint 
 of analysis 
 (or the \lq\lq{F-method}\rq\rq\
\cite{Eastwood60, KOSS, KP1}), 
 the multiplicity-two phenomenon
 in Observation \ref{obs:SL2} (2)
 can be derived from the fact
 that $\Omega_{\operatorname{sing}}$
 is of codimension two in $\Omega$
 and from the fact
 that
$\{
 \frac{\partial}{\partial \lambda_1}
 {\mathcal{RC}}_{\lambda_1,\lambda_2}^{\lambda_3}
, 
 \frac{\partial}{\partial \lambda_2}
 {\mathcal{RC}}_{\lambda_1,\lambda_2}^{\lambda_3}
\}
$
 forms a basis 
 in $H(\lambda_1,\lambda_2, \lambda_3)$
 when ${\mathcal{RC}}_{\lambda_1,\lambda_2}^{\lambda_3}=0$, 
 namely, 
 when $(\lambda_1,\lambda_2, \lambda_3)\in \Omega_{\operatorname{sing}}$.  
\item 
The basis given in Theorem \ref{thm:RC} (3) is different from 
 the basis in Remark \ref{rem:2.6} (1), 
 and clarifies the representation-theoretic reason
 for the multiplicity-two phenomenon
 as it is expressed
 as the composition
 of two intertwining operators.  
\item 
Theorem \ref{thm:RC} (3) implies a multiplicity-two phenomenon
 for Verma modules
 $M(\mu)=U({\mathfrak {g}}) \otimes_{U({\mathfrak {b}})} {\mathbb{C}}_{\mu}$
 for ${\mathfrak {g}}={\mathfrak {sl}}(2,{\mathbb{C}})$:
\[
 \dim_{\mathbb{C}} 
 \invHom{{\mathfrak {g}}}{M(-\lambda_3)}{M(-\lambda_1) \otimes M(-\lambda_2)}
=2
\quad
\text{for }
(\lambda_1,\lambda_2, \lambda_3)\in \Omega_{\operatorname{sing}}.  
\]
Again, 
 the tensor product 
 $M(-\lambda_1) \otimes M(-\lambda_2)$
 of Verma modules 
 decomposes into a multiplicity-free direct sum
 of irreducible ${\mathfrak {g}}$-modules
 for generic $\lambda_1, \lambda_2 \in {\mathbb{C}}$, 
 but not for singular parameters.  
See \cite[Part II]{KP1} for details.  
\item 
In turn, 
 we shall get a two-dimensional space
 of differential symmetry breaking operators
 at $\Omega_{\operatorname{sing}}$
 for principal series representations
 with respect to ${\mathrm{SL}}(2,{\mathbb{R}}) \times {\mathrm{SL}}(2,{\mathbb{R}}) 
\downarrow \operatorname{diag}({\mathrm{SL}}(2,{\mathbb{R}}))$, 
 see Remark \ref{rem:m2}
in Section \ref{sec:C}.  
\end{enumerate}
\end{rem}

\subsection{Symmetry breaking
 in conformal geometry}
\label{subsec:2.2}

In contrast
 to the localness property
 for symmetry breaking operators
 in the holomorphic setting 
 (Observation \ref{obs:SL2} (1)), 
 there exist non-local symmetry breaking operators
in a more general setting.  
We illustrate Stage C
 in the branching problem
 by an explicit construction
 and a complete classification 
 of all local and non-local symmetry breaking operators
 that arise from conformal geometry.  
In later sections,
 we explain a key idea
 of the proof
 (Section \ref{sec:C})
 and present potential settings
 where we could expect 
 that this example might serve
 as the prototype
 of analogous questions
 (Section \ref{sec:A2}).  
For full details
 of this subsection, 
 see the monograph \cite{xtsbon}
 joint with Speh.

For $\lambda \in {\mathbb{C}}$
 we denote by $I(\lambda)^{\infty}$
 the smooth (unnormalized) spherical principal series representation
 of $G={\mathrm{O}}(n+1,1)$.  
In our parametrization,
 $\lambda \in \frac n 2 + \sqrt{-1}{\mathbb{R}}$
 is the unitary axis, 
 $\lambda \in (0,n)$
 gives the complementary series representations,
 and $I(\lambda)^{\infty}$ contains 
 irreducible finite-dimensional representations
 as submodules
 for $\lambda \in \{0,-1,-2, \dots\}$
 and as quotients for $\lambda \in \{n ,n+1,n+2, \dots\}$.

We consider the restriction 
 of the representation $I(\lambda)^{\infty}$
 and its subquotients
 to the subgroup $G':={\mathrm{O}}(n,1)$.  
As we did for $I(\lambda)^{\infty}$, 
 we denote by $J(\nu)^{\infty}$
 for $\nu \in {\mathbb{C}}$, 
 the (unnormalized) spherical principal series representations
 of $G'={\mathrm{O}}(n,1)$.  
For $(\lambda,\nu) \in {\mathbb{C}}^2$, 
 we set 
\[
  H(\lambda, \nu):=\invHom{G'}{I(\lambda)^{\infty}}{J(\nu)^{\infty}}, 
\]
the space of (continuous) symmetry breaking operators.  
Similarly to Question \ref{que:SL2}, 
 we ask:
\begin{ques}
\label{que:On1}
\begin{enumerate}[{\rm (1)}]
\item
{\rm{(Stage B)}}\enspace
Find the dimension of $H(\lambda,\nu)$
 for $(\lambda,\nu)\in {\mathbb{C}}^2$.  
\item
{\rm{(Stage C)}}\enspace
Explicitly construct a basis for $H(\lambda,\nu)$.  
\item
{\rm{(Stage C)}}\enspace
Determine when $H(\lambda,\nu)$ contains 
 a differential operator.
\end{enumerate}

\end{ques}
The following is a complete answer to Question \ref{que:On1} (1).  
\begin{theorem}
\label{thm:On1}
\begin{enumerate}[{\rm (1)}]
\item 
For all $\lambda,\nu \in {\mathbb{C}}$, 
 we have $H(\lambda, \nu) \ne \{0\}$.  
\item 
$
\dim_{{\mathbb{C}}}
H(\lambda, \nu)
=
\begin{cases}
1
\qquad
& \text{if }\quad (\lambda,\nu) \in {\mathbb{C}}^2 \setminus L_{\operatorname{even}},   
\\
2
\qquad
& \text{if }\quad (\lambda,\nu) \in L_{\operatorname{even}},  
\end{cases}
$

where the \lq\lq{exceptional set}\rq\rq\
 $L_{\operatorname{even}}$
 is the discrete subset of ${\mathbb{C}}^2$
 defined by
\[
L_{\operatorname{even}}
:=
\{
(\lambda,\nu) \in {\mathbb{Z}}^2
:
 \lambda \le \nu \le 0,
\quad
 \lambda \equiv \nu \mod 2
\}.  
\]
\end{enumerate}
\end{theorem}

The role of $L_{\operatorname{even}}$
 in Theorem \ref{thm:On1}
 is similar to that of $\Omega_{\operatorname{sing}}$
 in Section \ref{subsec:SL2}.  
For Stage C, 
 we use the \lq\lq{$N$-picture}\rq\rq\
 of the principal series representations, 
 namely, 
 realize $I(\lambda)^{\infty}$ and $J(\nu)^{\infty}$
 in $C^{\infty}({\mathbb{R}}^n)$
 and $C^{\infty}({\mathbb{R}}^{n-1})$, 
 respectively.  
For $x \in{\mathbb{R}}^{n-1}$, 
 we set $|x|=(x_1^2+\cdots+x_{n-1}^2)^{\frac12}$.  
For $(\lambda, \nu) \in {\mathbb{C}}^2$
 satisfying $\operatorname{Re} (\nu - \lambda) \gg 0$
 and $\operatorname{Re}(\nu + \lambda) \gg 0$, 
 we construct explicitly
 a symmetry breaking operator
 (i.e.,  continuous $G'$-homomorphism)  from 
 $I(\lambda)^{\infty}$ to $J(\nu)^{\infty}$
 as an integral operator given by 
\begin{align}
\label{eqn:A}
  ({\mathbb{A}}_{\lambda,\nu} f)(y)
  :=
  &\int_{{\mathbb{R}}^{n}} |x_n|^{\lambda+\nu-n}
  (|x-y|^2 + x_n^2)^{-\nu}f(x,x_n)d x d x_n
\\
   =&\operatorname{rest}_{x_n=0}
    \circ
    (|x_n|^{\lambda+\nu-n}
(|x|^2+x_n^2)^{-\nu} 
\ast_{{\mathbb{R}^n}}
 f ).  
\notag
\end{align}

One might regard ${\mathbb{A}}_{\lambda,\nu}$ as a generalization
 of the Knapp--Stein intertwining operator
 ($G=G'$ case), 
 and also as the adjoint operator
 of a generalization
 of the Poisson transform.

The symmetry breaking operator
 ${\mathbb{A}}_{\lambda, \nu}$
 extends meromorphically with respect to the parameter 
 $(\lambda, \nu)$, 
 and if 
 we normalize ${\mathbb{A}}_{\lambda,\nu}$
 as
\[
  \widetilde{\mathbb{A}}_{\lambda,\nu}
  :=
\frac{1}{\Gamma(\frac{\lambda+\nu-n+1}{2})\Gamma(\frac{\lambda-\nu}{2})}
{\mathbb{A}}_{\lambda,\nu},   
\]
then $\widetilde {\mathbb{A}}_{\lambda,\nu}:I(\lambda)^{\infty} \to J(\nu)^{\infty}$
 is a continuous symmetry breaking operator
 that depends holomorphically
 on $(\lambda,\nu)$ in the entire complex plane ${\mathbb{C}}^2$, 
 and $\widetilde {\mathbb{A}}_{\lambda,\nu} \ne 0$
 if and only if $(\lambda,\nu)\not\in L_{\operatorname{even}}$
 (\cite[Theorem 1.5]{xtsbon}).

The singular set $L_{\operatorname{even}}$ is
 most interesting.  
To construct a symmetry breaking operator
 at $L_{\operatorname{even}}$, 
 we renormalize 
 $\widetilde {\mathbb{A}}_{\lambda,\nu}$
 for $\nu \in -{\mathbb{N}}$, 
 by 
\[
\widetilde{\widetilde {\mathbb{A}}}_{\lambda,\nu}
:=
\Gamma(\frac{\lambda-\nu}{2})
\widetilde {\mathbb{A}}_{\lambda,\nu}
=
\frac{1}{\Gamma(\frac{\lambda+\nu-n+1}{2})}
{\mathbb{A}}_{\lambda,\nu}.  
\]

In order to construct {\it{differential}} symmetry breaking operators, 
 we recall
 that the Gegenbauer polynomial 
 $C_l^{\alpha}(t)$
 for $l \in {\mathbb{N}}$ and $\alpha \in {\mathbb{C}}$
 is given by 
\[
   C_l^{\alpha}(t):=
  \sum_{k=0}^{[\frac l 2]}
  (-1)^k
  \frac{\Gamma(l-k+\alpha)}{\Gamma(\alpha) \Gamma(l-2k+1) k!}
  (2t)^{l-2k}.  
\]
We note that 
$C_l^{\alpha}(t) \equiv 0$
 if $l \ge 1$ and $\alpha =0,-1,-2, \dots, -[\frac{l-1}{2}]$.  
We renormalize 
 $C_l^{\alpha}(t)$
 by setting 
$
  \widetilde{C}_l^{\alpha}(t)
  :=
  \frac{\Gamma(\alpha)}{\Gamma(\alpha + [\frac{l+1}{2}])}
  C_l^{\alpha}(t), 
$
 so that $\widetilde{C}_l^{\alpha}(t)$ is a nonzero polynomial
 in $t$ of degree $l$
 for all $\alpha \in {\mathbb{C}}$
 and $l \in {\mathbb{N}}$.  
We inflate it 
 to a polynomial 
 of two variables $u$ and $v$
 by 
\[
  \widetilde C_{k}^{\alpha}(u,v)
  :=
  u^{\frac k 2} \widetilde C_{k}^{\alpha}(\frac{v}{\sqrt u}).  
\]
For instance,
 $\widetilde C_{0}^{\alpha}(u,v)=1$, 
 $\widetilde C_{1}^{\alpha}(u,v) = 2v$, 
 $\widetilde C_{2}^{\alpha}(u,v)=2(\alpha+1)v^2-u$, 
 etc. 
Substituting
 $u=-\Delta_{{\mathbb{R}^{n-1}}}=-\sum_{j=1}^{n-1}\frac{\partial^2}{\partial x_j^2}$
 and $v=\frac{\partial}{\partial x_n}$, 
 we get a differential operator
 of order $2l$:
\[
\widetilde {\mathbb{C}}_{\lambda,\nu}
:= 
\operatorname{rest}_{x_n=0} \circ
\,
 \widetilde C_{2l}^{\lambda -\frac{n-1}{2}}(-\Delta_{{\mathbb{R}^{n-1}}}, \frac{\partial}{\partial x_n}).  
\]

This closed formula
 of the differential operator $\widetilde{\mathbb{C}}_{\lambda, \nu}$
 was obtained
 by Juhl \cite{Juhl}
 (see also \cite{KOSS}
 for a short proof
 by the F-method, 
 and \cite{Eastwood60} for yet another proof
 by using the residue formula), 
 and the closed formula \eqref{eqn:A} of the symmetry breaking operator
 $\widetilde{\mathbb{A}}_{\lambda, \nu}$
 was obtained by Kobayashi and Speh \cite{xtsbon}.  

The following results answer 
 Question \ref{que:On1} (2) and (3); 
 see \cite[Theorems 1.8 and 1.9]{xtsbon}:
\begin{theorem}
\label{thm:Speh}
\begin{enumerate}[{\rm (1)}]
\item 
With notation as above, 
 we have
\[
H(\lambda,\nu)
=
\begin{cases}
{\mathbb{C}}\widetilde {\mathbb{A}}_{\lambda,\nu}
&\text{ if }(\lambda, \nu) \in {\mathbb{C}}^2 \setminus L_{\operatorname{even}} 
\\
{\mathbb{C}}\widetilde{\widetilde {\mathbb{A}}}_{\lambda,\nu}
\oplus 
{\mathbb{C}} \widetilde{\mathbb{C}}_{\lambda, \nu}
&\text{ if }(\lambda, \nu) \in L_{{\operatorname{even}}}.  
\end{cases}
\]
\item 
$H(\lambda,\nu)$ contains a nontrivial differential operator
 if and only if
 $\nu- \lambda =0,2,4,6, \dots$.  
In this case 
 $\widetilde {\mathbb{A}}_{\lambda,\nu}$ is proportional 
 to $\widetilde {\mathbb{C}}_{\lambda,\nu}$, 
 and the proportionality constant vanishes
 if and only if $(\lambda, \nu) \in L_{{\operatorname{even}}}$.  
\end{enumerate}
\end{theorem}
{}From Theorem \ref{thm:Speh} (2) and Theorem \ref{thm:On1} (1), 
 we have the following:

\begin{observation}
\begin{enumerate}[{\rm (1)}]
\item
Unlike the holomorphic setting 
 in Section \ref{subsec:SL2}, 
 the localness property fails.  
\item
Even if an irreducible smooth representation
 $\pi^{\infty}=I(\lambda)^{\infty}$ is unitarizable
 as a representation of $G$, 
 the condition 
$\invHom{G'}{\pi^{\infty}|_{G'}}{\tau^{\infty}} \ne \{0\}$
 does not imply 
 that the irreducible smooth representation
 $\tau^{\infty}=J(\nu)^{\infty}$
 is unitarizable
 as a representation of $G'$
 (see Section \ref{subsec:rep} for the terminology).  
\end{enumerate}
\end{observation}

For $\lambda \in \{n ,n+1,n+2, \dots\}$, 
 $I(\lambda)^{\infty}$ contains
 a unique proper infinite-dimensional closed $G$-submodule.  
We denote it by $A_{\mathfrak {q}}(\lambda-n)^{\infty}$, 
 which is 
 the Casselman--Wallach globalization 
 of Zuckerman's derived functor module
 $A_{\mathfrak {q}}(\lambda-n)$
 (see \cite{Vogan81, Vogan88})
 for some $\theta$-stable parabolic subalgebra
 ${\mathfrak {q}}$ of ${\mathfrak {g}}$.  
It is unitarizable
 (\cite{Vogan84, WaI})
 and has nonzero $({\mathfrak {g}}, K)$-cohomologies
 (Vogan--Zuckerman \cite{xvoza}).

By using the explicit formul{\ae}
 of symmetry breaking operators
 and certain identities 
 involving these operators,
 we can identify precisely the images
 of every subquotient
 of $I(\lambda)^{\infty}$
 under these operators.  
In particular,
 we obtain the following corollary 
 for the branching problem 
 of $A_{{\mathfrak{q}}}(\lambda)$ modules.  
We note
 that in this setting, 
 the restriction 
 $A_{{\mathfrak{q}}}(\lambda)|_{{\mathfrak{g}}'}$
 is not discretely decomposable
 as a $({\mathfrak {g}}', K')$-module 
 (Definition \ref{def:dd}).

\begin{corollary}
[{\cite[Theorem 1.2]{xtsbon}}]
With notation as above, 
 we have
\[
\dim_{\mathbb{C}}
\invHom{G'} {A_{\mathfrak {q}}(i)^{\infty}|_{G'}}{A_{\mathfrak {q}'}(j)^{\infty}}
=
\begin{cases}
1
\qquad
& \text{if }\quad i \ge j \quad \text{and} \quad i \equiv j \mod 2,  
\\
0
\qquad
& \text{if }\quad i < j \quad \text{and} \quad i \not\equiv j \mod 2.
\end{cases}
\]
\end{corollary}

There are some further applications
 of the explicit formul{\ae}
 \eqref{eqn:A}
 (Stage C in the branching problems).  
For instance, 
 J. M{\"o}llers and B. {\O}rsted
 recently found an interesting application 
 of the explicit formul\ae\ \eqref{eqn:A}
 to $L^p - L^q$ estimates
 of certain boundary-value problems, 
 and to some questions in automorphic forms
\cite{xmoor}.  

\section{Preliminary results and basic notation}
\label{sec:pre}
We review quickly some basic results
 on (infinite-dimensional) continuous representations
 of real reductive Lie groups
 and fix notation.  
There are no new results in this section.

\vskip 0.8pc
By a continuous representation $\pi$
 of a Lie group $G$
 on a topological vector space $V$
 we shall mean 
 that $\pi:G \to {\mathrm{GL}}_{\mathbb{C}}(V)$ is a group homomorphism from 
 $G$ into the group of invertible endomorphisms of $V$ 
 such that the induced map 
 $G \times V \to V$, 
 $(g,v) \mapsto \pi(g)v$
 is continuous.  
We say $\pi$ is a (continuous) Hilbert 
 [Banach, Fr{\'e}chet, $\cdots$]
 representation
 if $V$ is a Hilbert [Banach, Fr{\'e}chet, $\cdots$]
space.  
A continuous Hilbert representation $\pi$ of $G$
 is said to be a unitary representation
 when all the operators $\pi(g)$
 ($g \in G$)
 are unitary.

\subsection{Decomposition 
of unitary representations}
\label{sec:irreddeco}

One of the most distinguished features of \textit{unitary} representations 
is that they can be built up from the smallest objects, namely,
irreducible unitary representations.  
To be precise, 
let $G$ be a locally compact group.  
We denote by $\widehat G$
 the set of equivalence classes of irreducible
unitary representations of $G$
 (the {\it{unitary dual}}), 
 endowed with the Fell topology.

\begin{fact}
[Mautner--Teleman]
\label{fact:irreddeco}
For every unitary representation $\pi$ of a locally compact group $G$, 
there exist a Borel measure 
 $d \mu$ on $\dual G{}$
 and a measurable function 
 $n_{\pi} : \widehat{G} \to \mathbb{N} \cup \{ \infty \}$
such that $\pi$ is unitarily equivalent
 to the direct integral of irreducible unitary representations:

\begin{equation}
\label{eqn:3.1.1}
\pi \simeq \int_{\widehat{G}}^{\oplus} n_{\pi}(\sigma) \sigma
                  \, d\mu(\sigma), 
\end{equation}

where $n_\pi (\sigma) \sigma$ stands for the multiple of an
irreducible unitary representation $\sigma$ with multiplicity
$n_\pi(\sigma)$.
\end{fact}

The decomposition \eqref{eqn:3.1.1} is unique
 if $G$ is of type~I
in the sense of von Neumann algebras, 
in particular, 
 if $G$ (or $G'$ in later notation)
 is a real reductive Lie group
 or a nilpotent Lie group.  
Then the \textit{multiplicity function} $n_\pi$ is well-defined up to
a measure zero set with respect to $d\mu$.
We say
 that $\pi$ is \textit{multiplicity-free} if $n_\pi(\sigma)\le 1$
almost everywhere,
or equivalently,
if the ring of continuous $G$-endomorphisms of $\pi$ is commutative.

The decomposition \eqref{eqn:3.1.1} splits 
into a direct sum of the discrete and continuous parts:
\begin{equation}
\label{eqn:disco}
   \pi \simeq (\pi)_{\operatorname{disc}} \oplus (\pi)_{\operatorname{cont}}, 
\end{equation}
where $(\pi)_{\operatorname{disc}}$ is 
 a unitary representation defined
 on the maximal closed $G$-invariant subspace
 that is isomorphic to 
 a discrete Hilbert sum
 of irreducible unitary representations
and $(\pi)_{\operatorname{cont}}$ is its orthogonal complement.  

\begin{definition}
\label{def:discuni}
{\rm{
We say a unitary representation
 $\pi$ is {\it{discretely decomposable}}
 if $\pi=(\pi)_{\operatorname{disc}}$.  
}}
\end{definition}

\subsection{Continuous representations
 and smooth representations}
\label{subsec:rep}

We would like to treat non-unitary representations
 as well 
 for branching problems.  
For this we recall 
 some standard concepts
 of continuous representations of Lie groups.

Suppose $\pi$ is a continuous representation 
 of $G$ on a Banach space $V$.  
A vector $v \in V$ is said
 to be {\it{smooth}}
 if the map
 $G \to V$, 
 $g \mapsto \pi(g)v$ is of $C^{\infty}$-class.  
Let $V^{\infty}$ denote
 the space
 of smooth vectors
 of the representation $(\pi,V)$.  
Then $V^{\infty}$ carries 
 a Fr{\'e}chet topology
 with a family of semi-norms
$\|v\|_{i_1\cdots i_k}:=\|d\pi(X_{i_1}) \cdots d\pi(X_{i_k})v\|$, 
 where $\{X_1, \dots, X_n\}$ is a basis
 of the Lie algebra ${\mathfrak {g}}_0$ of $G$.  
Then $V^{\infty}$ is a $G$-invariant dense subspace
 of $V$, 
 and we obtain a continuous Fr{\'e}chet representation 
 $(\pi^{\infty}, V^{\infty})$
 of $G$.  
Similarly 
 we can define a representation $\pi^{\omega}$
 on the space $V^{\omega}$
 of analytic vectors.

Suppose now that $G$ is a real reductive linear Lie group, 
 $K$ a maximal compact subgroup of $G$,
 and ${\mathfrak {g}}$ the complexification
 of the Lie algebra ${\mathfrak {g}}_0$ of $G$.  
Let ${\mathcal{HC}}$ denote the category
 of Harish-Chandra modules
 whose objects and morphisms are $({\mathfrak {g}}, K)$-modules
 of finite length
 and $({\mathfrak {g}}, K)$-homomorphisms, 
 respectively.

Let $\pi$ be a continuous representation 
 of $G$ on a Fr{\'e}chet space $V$.  
Suppose 
 that $\pi$ is of finite length,
 namely,
 there are at most finitely many closed $G$-invariant subspaces
 in $V$.  
We say $\pi$ is {\it{admissible}}
 if 
\[
   \dim \operatorname{Hom}_K(\tau, \pi|_K)< \infty
\]
 for any irreducible finite-dimensional representation
 $\tau$ of $K$.  
We denote by $V_K$
 the space of $K$-finite vectors.  
Then $V_K \subset V^{\omega} \subset V^{\infty}$
 and the Lie algebra ${\mathfrak {g}}$
 leaves $V_K$ invariant.  
The resulting $({\mathfrak {g}}, K)$-module
 on $V_K$
 is called the underlying $({\mathfrak {g}}, K)$-module
 of $\pi$, 
 and will be denoted by $\pi_K$.

For any admissible representation $\pi$
 on a Banach space $V$, 
 the smooth representation $(\pi^{\infty}, V^{\infty})$
 depends only on the underlying 
 $({\mathfrak {g}}, K)$-module.  
We say $(\pi^{\infty}, V^{\infty})$
 is an {\it{admissible smooth representation}}.  
By the Casselman--Wallach globalization theory, 
 $(\pi^{\infty}, V^{\infty})$ has moderate growth, 
 and there is a canonical equivalence of categories
 between the category ${\mathcal{HC}}$
 of $({\mathfrak {g}}, K)$-modules
 of finite length 
 and the category of admissible smooth representations 
 of $G$ (\cite[Chapter 11]{WaI}).  
In particular,
 the Fr{\'e}chet representation
 $\pi^{\infty}$ is uniquely 
 determined by its underlying
 $({\mathfrak {g}}, K)$-module.  
We say $\pi^{\infty}$
 is the {\it{smooth globalization}}
 of $\pi_K \in {\mathcal{HC}}$.

For simplicity,
 by an {\it{irreducible smooth representation}}, 
 we shall mean an irreducible admissible smooth representation
 of $G$.  
We denote by $\dual {G}{smooth}$
 the set of equivalence classes
 of irreducible smooth representations 
 of $G$.  
Using the category ${\mathcal{HC}}$
 of $({\mathfrak {g}}, K)$-modules, 
 we may regard the unitary dual $\dual G{}$ 
as a subset of $\dual G {smooth}$.

\section{Two spaces:
 $\invHom{G'}{\tau}{\pi|_{G'}}$ and 
 $\invHom{G'}{\pi|_{G'}}{\tau}$
}
\label{sec:ep}
Given irreducible continuous representations $\pi$ of $G$
 and $\tau$ of a subgroup $G'$, 
 we may consider two settings for branching problems:
\par\noindent
{\bf{Case I.}}\enspace
(embedding)\enspace 
continuous $G'$-homomorphisms from $\tau$
 to $\pi|_{G'}$;
\par\noindent
{\bf{Case II.}}\enspace
(symmetry breaking)\enspace
 continuous $G'$-homomorphisms from $\pi|_{G'}$
 to $\tau$.

We write $\invHom{G'}{\tau}{\pi|_{G'}}$ and $\invHom{G'}{\pi|_{G'}}{\tau}$
 for the vector spaces
 of such continuous $G'$-homomorphisms, 
 respectively.  
Needless to say, 
 the existence
 of such $G'$-intertwining operators
 depends on the topology of the representation spaces
 of $\pi$ and $\tau$.

Cases I and II are related
 to each other 
 by taking contragredient representations:
\begin{align*}
\invHom{G'}{\tau}{\pi|_{G'}} \subset 
& \invHom{G'}{\pi^{\vee}|_{G'}}{\tau^{\vee}}, 
\\
\invHom{G'}{\pi|_{G'}}{\tau} 
\subset & 
\invHom{G'}{\tau^{\vee}}{\pi^{\vee}|_{G'}}.  
\end{align*}
Thus they are equivalent
 in the category
 of unitary representations
 (see Theorem \ref{thm:taupi} (3)).  
Furthermore, 
 we shall use a variant of the above duality
 in analyzing differential symmetry breaking operators (Case II)
by means of \lq\lq discretely decomposable restrictions\rq\rq\
of Verma modules (Case I);
 see the duality \eqref{eqn:Vdual} for the proof 
 of Theorem \ref{thm:finVerma} below.

On the other hand,
 it turns out
that Cases I and II are significantly different
 if we confine ourselves
 to irreducible smooth representations
 (see Section \ref{subsec:rep}).  
Such a difference 
 also arises in an analogous problem
 in the category ${\mathcal{H C}}$
 of Harish-Chandra modules
 where  no topology is specified.

Accordingly, 
 we shall discuss some details for Cases I and II separately,
 in Sections \ref{sec:A} and \ref{sec:A2}, 
 respectively. 

\subsection{$K$-finite vectors and $K'$-finite vectors}
\label{subsec:gkemb}
Let $G$ be a real reductive linear Lie group,
 and $G'$ a reductive subgroup.  
We take maximal compact subgroups $K$ and $K'$
 of $G$ and $G'$, 
 respectively,
 such that $K' = K \cap G'$.

We recall that for an admissible representation $\pi$ of $G$
 on a Banach space $V$, 
 any $K$-finite vector is contained
 in $V^{\infty}$, 
and the underlying $({\mathfrak {g}}, K)$-module $\pi_K$ is defined on 
\[
   V_K := V_{\text{$K$-finite}} 
   \qquad
   (\subset V^{\infty}).  
\]
When we regard $(\pi, V)$
 as a representation of the subgroup $G'$
 by restriction,
 we denote by $(V|_{G'})^{\infty}$ the space
 of smooth vectors
 with respect to the $G'$-action, 
 and write $(\pi|_{G'})^{\infty}$ for the continuous representation of $G'$
 on $(V|_{G'})^{\infty}$.  
In contrast to the case $G=G'$, 
 we remark
 that $K'$-finite vectors are not necessarily
 contained in $(V|_{G'})^{\infty}$
 if $G' \subsetneqq G$, 
 because the $G'$-module $(\pi|_{G'}, V|_{G'})$ 
 is usually not of finite length.  
Instead, 
 we can define a $({\mathfrak {g}}', K')$-module
 on 
\[
  V_{K'} := V_{\text{$K'$-finite}}
                        \cap  
                        (V|_{G'})^{\infty}, 
\]
which we denote simply by $\pi_{K'}$.  
Obviously we have the following inclusion relations:
\begin{alignat}{3}
& V_K && \subset \,\,&& \,\,\,\,\,V_{K'}
\notag
\\
&\, \cap &&        && \,\, \,\,\,\,\cap
\label{eqn:VKinfty}
\\
& V^{\infty} && \subset \,\,&& (V|_{G'})^{\infty} \subset V
\notag
\end{alignat}

None of them coincides in general 
 (e.g., 
 $V_K=V_{K'}$
 if and only if $\pi_K$ is discretely decomposable
 as $({\mathfrak {g}}', K')$-module, 
as we shall see in Theorem \ref{thm:deco98} below.

We set
\begin{align*}
 H_K(\tau, \pi):=&\invHom{{\mathfrak {g}}', K'}{\tau_{K'}}{\pi_K|_{\mathfrak {g}'}}, 
\\
 H_{K'}(\tau, \pi):=&
 \invHom{{\mathfrak {g}}', K'}{\tau_{K'}}{\pi_{K'}|_{\mathfrak {g}'}}.  
\end{align*}
According to the inclusion relation
 \eqref{eqn:VKinfty}, 
 for irreducible representations $\tau$ of $G'$ we have:
\[
   H_K(\tau, \pi) \subset H_{K'}(\tau, \pi).  
\]
In the case where $\pi$ is a unitary representation of $G$, 
 the latter  captures
discrete summands in the branching law
 of the restriction $\pi|_{G'}$
 (see, Theorem \ref{thm:taupi} (3)), 
 whereas the former vanishes
 even if the latter is nonzero
 when the continuous part
 $(\pi|_{G'})_{\operatorname{cont}}$ is not empty
 (see Theorem \ref{thm:deco98}).  
The spaces
 of continuous $G'$-homomorphisms
 such as $\invHom{G'}{\tau}{\pi|_{G'}}$
 or $\invHom{G'}{\tau^{\infty}}{\pi^{\infty}|_{G'}}$
 are in between.

\vskip 1pc
We begin with a general result:
\begin{theorem}
\label{thm:taupi}
Suppose that $\pi$ and $\tau$
 are admissible irreducible Banach representations
 of $G$ and $G'$.  
\begin{enumerate}[{\rm (1)}]
\item 
We have natural inclusions and an isomorphism:
\begin{multline}
H_K(\tau, \pi)
\hookrightarrow
\invHom{G'}{\tau^{\infty}}{\pi^{\infty}|_{G'}}
\\
\hookrightarrow
\invHom{G'}{\tau^{\infty}}{(\pi|_{G'})^{\infty}}
\overset \sim \to 
H_{K'}(\tau, \pi).  
\label{eqn:taupi}
\end{multline}
\item 
 There are canonical injective homomorphisms:
\begin{multline}
\invHom {G'}{\pi|_{G'}}{\tau}
 \hookrightarrow 
\invHom{G'}{\pi^{\infty}|_{G'}}{\tau^{\infty}}
\\
 \hookrightarrow 
 \invHom{G'}{\pi^{\omega}|_{G'}}{\tau^{\omega}}
\hookrightarrow
 \invHom{{\mathfrak {g}}', K'}{\pi_K}{\tau_{K'}}.  
\label{eqn:Hsmooth}
\end{multline}
\item 
{\rm{(unitary case)}}
If $\tau$ and $\pi$ are irreducible unitary representations
 of $G'$ and $G$, 
respectively,
 then we have natural isomorphisms
 (where the last isomorphism is conjugate linear):
\begin{multline}
\label{eqn:unid}
H_{K'}(\tau, \pi)
\overset \sim \leftarrow
\invHom{G'}{\tau^{\infty}}{(\pi|_{G'})^{\infty}}
\\
\overset \sim \leftarrow 
\invHom{G'}{\tau}{\pi|_{G'}}
\simeq
\invHom{G'}{\pi|_{G'}}{\tau}.  
\end{multline}
We write $m_{\pi}(\tau)$
 for the dimension of one of (therefore, any of)
 the terms in \eqref{eqn:unid}.  
Then the discrete part
 of the restriction $\pi|_{G'}$
 (see Definition \ref{def:discuni})
 decomposes discretely as
\[
  (\pi|_{G'})_{\operatorname{disc}}
  \simeq 
  {\sum_{\tau \in \widehat {G'}}}{}^{\oplus}  
  m_{\pi}(\tau)\tau.   
\]
\end{enumerate}
\end{theorem}

\begin{rem}
Even if $\pi$ and $\tau$ are irreducible unitary representations
 of $G$ and $G'$, 
 respectively,
 the canonical injective homomorphism
\begin{equation}
\label{eqn:pisurj}
\invHom{G'}{\pi|_{G'}}{\tau} \hookrightarrow
\invHom{G'}{\pi^{\infty}|_{G'}}{\tau^{\infty}}
\end{equation}
is not surjective in general.

In fact,
 we can give an example
 where the canonical homomorphism \eqref{eqn:pisurj}
 is not surjective
 by using the classification 
 of $\invHom{G'}{\pi^{\infty}|_{G'}}{\tau^{\infty}}$
 for the pair $(G,G')=(O(n+1,1),O(n,1))$
 in Section \ref{subsec:2.2}
 as follows.  
Recall  
 $\invHom{G'}{I(\lambda)^{\infty}|_{G'}}{J(\nu)^{\infty}}$ $\ne \{0\}$
 for all $(\lambda,\nu) \in {\mathbb{C}}^2$
 with the notation therein.  
However,
 for a fixed $\pi \in \dual G{}$, 
 there exist 
 at most countably many $\tau \in \dual{G'}{}$
 that occur
 in the discrete part
 of the restriction $\pi|_{G'}$, 
 and therefore 
$
   \{\tau \in \dual{G'}{}: \invHom{G'}{\pi|_{G'}}{\tau} \ne \{0\}\}
$
 is an infinite set because we have the following bijection:  
\[
\{
   \tau \in \dual{G'}{}
:
  \invHom{G'}{\pi|_{G'}}{\tau} \ne \{0\}
\}
\simeq
\{
   \tau \in \dual{G'}{}
:
  \invHom{G'}{\tau}{\pi|_{G'}} \ne \{0\}
\}.  
\]
Hence, 
by taking $\pi^{\infty} = I(\lambda)^{\infty}$
 for a fixed $\lambda \in \frac n 2+ \sqrt{-1} {\mathbb{R}}$
 (unitary axis)
 or $\lambda \in (0,n)$
 (complementary series), 
 we see
 that the canonical homomorphism
 \eqref{eqn:pisurj}
 must be zero
 when we take $\tau^{\infty}$
 to be a representation $I(\nu)^{\infty}$
 for $\nu \in {\mathbb{C}}$
 such that $\nu \not \in \frac{n-1}{2}+ \sqrt{-1} {\mathbb{R}}$
 and $\nu \not \in {\mathbb{R}}$.  
\end{rem}

Let us give a proof of Theorem \ref{thm:taupi}.  
\begin{qedproof}
(1)\enspace
To see the first inclusion,
 we prove
 that any $({\mathfrak {g}}', K')$-homomorphism
 $\iota:\tau_{K'} \to \pi_K|_{{\mathfrak {g}}'}$
 extends to a continuous map
 $\tau^{\infty} \to \pi^{\infty}|_{G'}$.  
We may assume
 that $\iota$ is nonzero,
 and therefore,
 is injective.  
Since $\iota(\tau_{K'}) \subset \pi_K \subset \pi^{\infty}$, 
 we can define a Fr{\'e}chet space $W$
 to be the closure of $\iota(\tau_{K'})$ in $\pi^{\infty}$, 
 on which $G'$ acts continuously.  
Its underlying $({\mathfrak {g}}', K')$-module
 is isomorphic to $\iota(\tau_{K'}) \simeq \tau_{K'}$.

Since the continuous representation $\pi^{\infty}$ of $G$ 
 is of moderate growth,
 the  Fr{\'e}chet representation $W$ of the subgroup $G'$
 is also of moderate growth.  
By the Casselman--Wallach globalization theory,
 there is a $G'$-homomorphism $\tau^{\infty} \overset \sim \to 
\overline{\iota(\tau_{K'})}$ $(=W)$
 extending the $({\mathfrak {g}}', K')$-isomorphism
 $\iota:\tau_{K'} \overset \sim\to \iota(\tau_{K'})$.  
Hence we have obtained a natural map 
 $\invHom{{\mathfrak {g}}', K'}{\tau_{K'}}{\pi_K|_{{\mathfrak {g}}'}}
 \to \invHom{G'}{\tau^{\infty}}{\pi^{\infty}|_{G'}}$, 
 which is clearly injective 
 because $\tau_{K'}$ is dense in $\tau^{\infty}$.

The second inclusion is obvious.

To see the third inclusion,
 it suffices to show
that any $\iota \in \invHom{{\mathfrak {g}}', K'}{\tau_{K'}}{\pi_{K'}|_{{\mathfrak {g}}'}}$
 extends to a continuous $G'$-homomorphism from 
$\tau^{\infty}$ to $(\pi|_{G'})^{\infty}$.  
Since $\tau_{K'}$ is an irreducible $({\mathfrak {g}}', K')$-module, 
 $\iota$ is injective
 unless $\iota$ is zero
 and $\iota(\tau_{K'})$ is isomorphic to $\tau_{K'}$
 as $({\mathfrak {g}}', K')$-modules.

Let $V$ be the Banach space
 on which $G$ acts via $\pi$, 
 and $W_1$ and $W_2$ the closures of $\iota(\tau_{K'})$
 in the Banach space $V$
 and the Fr{\'e}chet space $(V|_{G'})^{\infty}$, 
 respectively.  
Then the underlying $({\mathfrak {g}}', K')$-modules
 of $W_1$ and $W_2$ are both isomorphic to $\tau_{K'}$.  
Moreover,
 $W_2 \subset W_1 \cap (V|_{G'})^{\infty}$ by definition, 
 and $W_2$ is closed in $W_1 \cap (V|_{G'})^{\infty}$
 with respect to the Fr{\'e}chet topology. 
Since the subspace $\iota(\tau_{K'})$ of $W_2$ is dense
 in $W_1 \cap (V|_{G'})^{\infty}$, 
 we conclude
 that $W_2$ coincides with $W_1 \cap (V|_{G'})^{\infty}$, 
 which is the Casselman--Wallach globalization
 of the $({\mathfrak {g}}', K')$-module
 $\iota(\tau_{K'})\simeq \tau_{K'}$.  
By the uniqueness
 of the Casselman--Wallach globalization \cite[Chapter 11]{WaI}, 
the $({\mathfrak {g}}', K')$-isomorphism
 $\tau_{K'} \overset {\sim} \to \iota(\tau_{K'})$
 extends to an isomorphism
 between Fr{\'e}chet $G'$-modules 
 $\tau^{\infty} \overset {\sim} \to W_2
 (=W_1 \cap (V|_{G'})^{\infty})$.  
\par\noindent
(2)\enspace
If $\iota:\pi|_{G'} \to \tau$ is a continuous $G'$-homomorphism,
 then 
\[
         \iota(\pi^{\infty}|_{G'}) 
\subset  \iota((\pi|_{G'})^{\infty})
\subset  \tau^{\infty}, 
\]
 and thus we have obtained a continuous $G'$-homomorphism
 $\iota^{\infty}:\pi^{\infty}|_{G'} \to \tau^{\infty}$
 between Fr{\'e}chet representations.  
Furthermore $\iota \mapsto \iota^{\infty}$ is injective
 because $V^{\infty}$ is dense in $V$.  
This shows the first inclusive relation
 of the statement (2).  
The proof for other inclusions are similar.  
\par\noindent
(3)\enspace
The last isomorphism in \eqref{eqn:unid} is given by taking the adjoint operator.  
The other isomorphisms are easy to see. 
The last statement follows from the fact
that if $\varphi \in \invHom{G'}{\tau}{\pi|_{G'}}$
 then $\varphi$ is a scalar multiple
 of an {\it{isometric}} $G'$-homomorphism.  
\end{qedproof}

The terms in \eqref{eqn:taupi}
 do not coincide in general.  
In order to clarify
 when they coincide, 
 we recall from \cite{xkInvent98}
 the notion of discrete decomposability 
 of ${\mathfrak {g}}$-modules.

\begin{definition}
\label{def:dd}
{\rm 
A $(\mathfrak{g},K)$-module $X$
 is said to be {\it{discretely decomposable as
 a $(\mathfrak{g}',K')$-module}}
 if there is a filtration $\{X_i\}_{i\in \mathbb{N}}$ of
 $(\mathfrak{g}',K')$-modules such that 
\begin{itemize}
\item $\bigcup_{i\in \mathbb{N}} X_i=X$ and
\item $X_i$ is of finite length as a $(\mathfrak{g}',K')$-module
 for any $i\in \mathbb{N}$.
\end{itemize}
The idea was
 to exclude \lq\lq{hidden continuous spectrum}\rq\rq\
 in an algebraic setting,
 and discrete decomposability here
 does not imply complete reducibility.  
Discrete decomposability is preserved
 by taking submodules,
 quotients,
 and the tensor product 
 with finite-dimensional representations.  
}
\end{definition}

\begin{rem}[see {\cite[Lemma 1.3]{xkInvent98}}]
{\rm 
Suppose that $X$ is a unitarizable $(\mathfrak{g},K)$-module.
Then $X$ is discretely decomposable  as a $(\mathfrak{g}',K')$-module
if and only if 
$X$ is isomorphic to 
an algebraic direct sum of irreducible $(\mathfrak{g}',K')$-modules.
}
\end{rem}

We get much stronger results
 than Theorem \ref{thm:taupi}
in this setting:
\begin{theorem}
[discretely decomposable case]
\label{thm:deco98}
Assume $\pi$ is an irreducible admissible representation of $G$
 on a Banach space $V$.  
Let $\pi_K$ be the underlying $({\mathfrak {g}}, K)$-module.  
Then the following five conditions
 on the triple $(G,G',\pi)$ are equivalent:
\begin{enumerate}[{\rm (i)}]
\item
There exists at least one irreducible $({\mathfrak {g}}', K')$-module $\tau_{K'}$
 such that 

$
   \invHom{{\mathfrak {g}}', K'}{\tau_{K'}}{\pi_K}\ne \{0\}.  
$
\item
$\pi_K$ is discretely decomposable as a $({\mathfrak {g}}', K')$-module
 (see Definition \ref{def:dd}).  
\item
All the terms in \eqref{eqn:taupi} are the same
 for any irreducible admissible Banach representation $\tau$ of $G'$.  
\item
All the terms in \eqref{eqn:taupi} are 
 the same for some irreducible admissible Banach representation $\tau$ of $G'$. \item
$V_K = V_{K'}$.   
\end{enumerate}
Moreover, 
if $(\pi,V)$ is a unitary representation,
 then one of (therefore, any of) the equivalent conditions
 (i) -- (v) implies
 that the continuous part $(\pi|_{G'})_{\operatorname{cont}}$
 of the restriction $\pi|_{G'}$
 is empty.  
\end{theorem}
\begin{qedproof}
See \cite{xkInvent98}
 for the first statement, 
 and \cite{xkdecoaspm} for the second statement.  
\end{qedproof}

\subsection{Some observations on 
 $\invHom{G'}{\tau^{\infty}}{\pi^{\infty}|_{G'}}$ and $\invHom{G'}{\pi^{\infty}|_{G'}}{\tau^{\infty}}$}
\label{subsec:2}

For a unitary representation $(\pi,V)$
 of $G$,
 Fact \ref{fact:irreddeco} gives an irreducible decomposition 
 of the restriction $\pi|_{G'}$ 
 into irreducible {\it{unitary}} representations of $G'$.  
However, symmetry breaking operators may exist
 between unitary and non-unitary
 representations:
\begin{observation}
\label{obs:unitary}
Suppose $\pi$ is a unitary representation 
 of $G$, 
 and $(\tau, W)$ an irreducible admissible representation
 of a reductive subgroup $G'$.  
\begin{enumerate}[{\rm (1)}]
\item 
If $\invHom{G'}{\tau^{\infty}}{\pi^{\infty}|_{G'}} \ne \{0\}$, 
 then $\tau^{\infty}$ is unitarizable.  
Actually, 
 $\tau$ occurs as a discrete part
 of $(\pi|_{G'})_{\operatorname{disc}}$
 (see \eqref{eqn:disco}).  
\item 
It may well happen 
 that 
$\invHom{G'}{\pi^{\infty}|_{G'}}{\tau^{\infty}} \ne \{0\}$
even when $\tau^{\infty}$ is not unitarizable.  
\end{enumerate}
\end{observation}
In fact, 
the first assertion is obtained
 by taking the completion of $\varphi(W^{\infty})$
 in the Hilbert space $V$ for 
 $\varphi \in \invHom{G'}{\tau^{\infty}}{\pi^{\infty}|_{G'}}$
 as in the proof of Theorem \ref{thm:taupi} (3), 
 where we considered the case $(\pi|_{G'})^{\infty}$
 instead of $\pi^{\infty} |_{G'}$.  
Theorem \ref{thm:Speh} gives an example
 of Observation \ref{obs:unitary} (2).

Here is another example
 that indicates a large difference
 between the two spaces, 
 $\invHom{G'}{\tau^{\infty}}{\pi^{\infty}|_{G'}}$
 and $\invHom{G'}{\pi^{\infty}|_{G'}}{\tau^{\infty}}$.  
\begin{examp}
{\rm{
Suppose $G$ is a real simple connected Lie group, 
 and $G'$ is a noncompact closed subgroup of $G$.  
Let $\pi$ be any irreducible unitary representation
 such that $\dim \pi= \infty$
 and $\invHom{G}{\pi^{\infty}}{C^{\infty}(G/G')}\ne \{0\}$.  
Then 
 by Howe--Moore \cite{xhowemoore}
 we have
\[
\invHom{G'}{{\bf{1}}}{\pi^{\infty}|_{G'}}
=\{0\}
\ne 
\invHom{G'}{\pi^{\infty}|_{G'}}{{\bf{1}}}.  
\]
}}
\end{examp}

\section{Features of the restriction, I
 : $\invHom{G'}{\tau}{\pi|_{G'}}$ (embedding)}
\label{sec:A}
In this section,
 we discuss Case I in Section \ref{sec:ep},
 namely
 $G'$-homomorphisms from irreducible $G'$-modules $\tau$
 into irreducible $G$-modules $\pi$.  
We put emphasis on its algebraic analogue
 in the category ${\mathcal{HC}}$ of Harish-Chandra modules.

The goals of this section are
\begin{enumerate}[(1)]
\item 
(criterion)\enspace
to find a criterion 
 for the triple $(G,G',\pi)$
 such that 
\begin{equation}
\label{eqn:120}
   \invHom{{\mathfrak{g}}',K'}{\tau_{K'}}{\pi_{K}|_{\mathfrak {g}'}} \ne \{0\} 
\quad
\text{for some $\tau$};
\end{equation}
\item 
(classification theory) \enspace
to classify the pairs $(G,G')$
 of reductive groups 
 for which \eqref{eqn:120} occurs
 for at least one infinite-dimensional $\pi \in \dual G{}$.  
\end{enumerate}
We also discuss recent progress
 in this direction
 as a refinement of (2):
\begin{enumerate}[{\rm (2)}$'$]
\item (classification theory)\enspace
Classify the triples $(G,G',\pi)$
 for which \eqref{eqn:120} occurs
 in typical cases
 (e.g., 
 $\pi_K$ is Zuckerman's $A_{\mathfrak {q}}(\lambda)$ module, 
 or a minimal representation).  
\end{enumerate}

In Section \ref{sec:C}
 we shall explain two new applications
 of discretely decomposable restrictions:
one is a dimension estimate
 of differential symmetry breaking operators
 (Theorem \ref{thm:finVerma}), 
 and the other is a proof
 of the \lq\lq{localness property}\rq\rq\
 of symmetry breaking operators
 (Theorem \ref{thm:local});
 see Observation \ref{obs:SL2} (1).

\subsection{Criteria for discrete decomposability 
 of restriction}
\label{subsec:crideco}
We review a necessary and sufficient condition
 for the restriction  of Harish-Chandra modules
 to be discretely decomposable
 (Definition \ref{def:dd}),
 which was established in
 \cite{xkAnn98} and \cite{xkInvent98}.

An associated variety ${\mathcal{V}}_{\mathfrak{g}}(X)$
 is a coarse approximation of 
 the ${\mathfrak {g}}$-modules $X$, 
 which we recall now from Vogan \cite{xvoganass}.  
We shall use the associated variety 
 for the study of the restrictions
 of Harish-Chandra modules.

Let $\{U_j(\mathfrak{g})\}_{j \in {\mathbb{N}}}$
 be the standard increasing filtration
 of the universal enveloping algebra
 $U(\mathfrak{g})$.
Suppose $X$ is a finitely generated $\mathfrak{g}$-module.
A filtration $\bigcup_{i\in \mathbb{N}} X_i=X$
 is called a {\it good filtration} if
 it satisfies the following conditions:
\begin{itemize}
\item $X_i$ is finite-dimensional for any $i\in \mathbb{N}$;
\item $U_j(\mathfrak{g})X_i\subset X_{i+j}$ for any $i,j\in\mathbb{N}$;
\item There exists $n$ such that $U_j(\mathfrak{g})X_i=X_{i+j}$
 for any $i\geq n$ and $j\in \mathbb{N}$.
\end{itemize}

The graded algebra
 $\operatorname{gr} U(\mathfrak{g}):=
 \bigoplus_{j\in\mathbb{N}}U_j(\mathfrak{g})/U_{j-1}(\mathfrak{g})$
is isomorphic to the symmetric algebra $S(\mathfrak{g})$
 by the Poincar\'{e}--Birkhoff--Witt theorem and we 
 regard the graded module $\operatorname{gr}X:=\bigoplus_{i\in\mathbb{N}} X_i/X_{i-1}$
 as an $S(\mathfrak{g})$-module.  
Define
\begin{align*}
\operatorname{Ann}_{S(\mathfrak{g})}(\operatorname{gr}X):=
\{f\in S(\mathfrak{g}) : fv=0 \text{\ for any\ } v\in\operatorname{gr}X\},
\\
{\cal V}_{\mathfrak{g}}(X):=
\{x\in \mathfrak{g}^* : f(x)=0 \text{\ for any\ }
 f\in\operatorname{Ann}_{S(\mathfrak{g})}(\operatorname{gr}X)\}.
\end{align*}
Then ${\cal V}_{\mathfrak{g}}(X)$ does not depend on
 the choice of good filtration and is called
 the {\it associated variety} of $X$.
We denote by ${\mathcal{N}}({\mathfrak {g}}^{\ast})$
 the nilpotent variety 
 of the dual space ${\mathfrak {g}}^{\ast}$.  
We have then the following basic properties
 of the associated variety \cite{xvoganass}.

\begin{lemma}
\label{lem:ass}
Let $X$ be a finitely generated $\mathfrak{g}$-module.
\begin{enumerate}[{\rm (1)}]
\item 
If $X$ is of finite length, then
 ${\cal V}_{\mathfrak{g}}(X)\subset{\cal N}(\mathfrak{g}^*)$.
\item 
 ${\cal V}_{\mathfrak{g}}(X)=\{0\}$ if and only if $X$ is finite-dimensional.
\item 
Let $\mathfrak{h}$ be a Lie subalgebra of $\mathfrak{g}$.
Then ${\cal V}_{\mathfrak{g}}(X)\subset \mathfrak{h}^\perp$
 if $\mathfrak{h}$ acts locally finitely on $X$,
 where 
 $\mathfrak{h}^\perp:=
 \{x \in {\mathfrak {g}}^{\ast}:
 x|_{\mathfrak{h}}=0\}$.
\end{enumerate}
\end{lemma}

(1) and (3) imply that
if $X$ is a $(\mathfrak{g},K)$-module of finite length,
 then ${\cal V}_\mathfrak{g}(X)$ is
 a $K_\mathbb{C}$-stable closed subvariety of
 ${\cal N}(\mathfrak{p}^*)$
 because $\mathfrak{k}^\perp=\mathfrak{p}^{\ast}$.

Dual to the inclusion ${\mathfrak{g}}' \subset {\mathfrak{g}}$
 of the Lie algebras, 
 we write
\[\operatorname{pr}:\mathfrak{g}^* \to (\mathfrak{g}')^* \]
 for the restriction map.

One might guess 
 that irreducible summands
 of the restriction $\pi|_{G'}$ would be
 \lq\lq{large}\rq\rq\
 if the irreducible representation $\pi$ of $G$ is \lq\lq{large}\rq\rq.  
The following theorem shows
 that such a statement holds
 if the restriction
 of the Harish-Chandra module
 is discretely decomposable
 (Definition \ref{def:dd}), 
however, 
 it is false in general 
 (see Counterexample \ref{ex:small} below).  
\begin{fact}
\label{nec}
Let $X$ be an irreducible $(\mathfrak{g},K)$-module.   
\begin{enumerate}[{\rm (1)}]
\item 
If $Y$ is an irreducible $(\mathfrak{g}',K')$-module
such that 
 $\invHom{\mathfrak{g}',K'}{Y}{X|_{\mathfrak {g}'}} \ne \{0\}$, 
then 
\[
   \operatorname{pr}({\mathcal {V}}_{\mathfrak {g}}(X)) \subset {\mathcal {V}}_{\mathfrak {g}'}(Y).  
\]
\item 
If 
$Y^{(j)}$ are irreducible $(\mathfrak{g}',K')$-modules
 such that 
$\invHom{\mathfrak{g}',K'}{Y^{(j)}}{X|_{\mathfrak {g}'}} \ne \{0\}$
 $(j=1,2)$, 
then 
\[
   {\mathcal {V}}_{\mathfrak {g}'}(Y_1) = {\mathcal {V}}_{\mathfrak {g}'}(Y_2).  
\]
In particular,
 the Gelfand--Kirillov dimension $\operatorname{GK-dim}(Y)$ 
of all irreducible $(\mathfrak{g}',K')$-submodules $Y$ of 
 $X|_{\mathfrak {g}'}$ are the same.  
\item 
{\rm{(necessary condition {\cite[Corollary 3.5]{xkInvent98}})}}\enspace
If $X$ is discretely decomposable
 as a $(\mathfrak{g}',K')$-module, 
 then 
$\operatorname{pr} ({\cal V}_\mathfrak{g} (X))\subset
 {\cal N}((\mathfrak{g}')^*)$,
 where ${\cal N}((\mathfrak{g}')^*)$ is the nilpotent variety
 of ${(\mathfrak{g}')}^*$.
\end{enumerate}
\end{fact}

An analogous statement fails
 if we replace
 $\invHom{{\mathfrak{g}}', K'}{\tau_{K'}}{\pi_K|_{{\mathfrak{g}}'}}$
 by the space 
 $\invHom{G'}{\tau}{\pi|_{G'}}$
 of continuous $G'$-intertwining operators:
\begin{false}
\label{false:prv}
Let $\pi$ be an irreducible unitary representation
 of a real reductive Lie group $G$.  
\begin{enumerate}[{\rm (1)}]
\item 
If $\tau \in \dual{G'}{}$ satisfies $\invHom{G'}{\tau}{\pi|_{G'}} \ne \{0\}$, 
 then 
$
   \operatorname{pr}({\mathcal {V}}_{\mathfrak {g}}(\pi_K))
 \subset {\mathcal {V}}_{\mathfrak {g}'}(\tau_{K'}).  
$
\item 
If $\tau^{(j)}$ $\in \dual{G'}{}$
 satisfy $\invHom{G'}{\tau^{(j)}}{\pi|_{G'}} \ne \{0\}$ ($j=1,2$), 
 then 
$
   {\mathcal {V}}_{\mathfrak {g}'}(\tau_{K'}^{(1)})
 = {\mathcal {V}}_{\mathfrak {g}'}(\tau_{K'}^{(2)}).  
$
\end{enumerate}
\end{false}

Here are counterexamples 
 to the \lq\lq{False Statement \ref{false:prv}}\rq\rq:
\begin{counterexample}
\label{ex:small}
\begin{enumerate}[{\rm (1)}]
\item 
There are many triples
 $(G,G',\pi)$
 such that $\pi \in \dual G{}$
 satisfies $(\pi|_{G'})_{\operatorname{cont}} \ne 0$;
 see \cite[Introduction]{xkInvent94}, 
 \cite[Section 3.3]{deco-euro}, 
 and Theorem \ref{thm:tensorP}, 
 for instance.  
In this case, 
$
\operatorname{pr}({\mathcal {V}}_{\mathfrak {g}}(\pi_K))
 \not\subset {\mathcal {V}}_{\mathfrak {g}'}(\tau_{K'})
$
 for any $\tau \in \dual{G'}{}$
 by Fact \ref{nec} (3).  
\item 
Let $(G, G')=(G_1 \times G_1, \operatorname{diag}(G_1))$
with $G_1={\mathrm{Sp}}(n,{\mathbb{R}})$ ($n \ge 2$).  
Take an irreducible unitary spherical principal series representation $\pi_1$
 induced from the Siegel parabolic subgroup of $G_1$, 
 and set $\pi = \pi_1 \boxtimes \pi_1$.  
Then there exist discrete series representations $\tau^{(1)}$ and $\tau^{(2)}$
 of $G'$ $(\simeq {\mathrm{Sp}}(n,{\mathbb{R}}))$, 
 where $\tau^{(1)}$ is a holomorphic discrete series representation
 and $\tau^{(2)}$ is a non-holomorphic discrete series representation, 
such that 
\[
  \invHom{G'}{\tau^{(j)}}{\pi} \ne \{0\}
  \quad
 (j=1,2)
\quad
\text{and}
\quad
\operatorname{GK-dim}(\tau^{(1)})
<
\operatorname{GK-dim}(\tau^{(2)}).  
\]
\end{enumerate}
\end{counterexample}

In fact, 
 it follows from Theorem \ref{thm:tensorP} below
 that $\invHom{G'}{\tau}{\pi} \ne \{0\}$
 if and only if $\tau$ is a discrete series representation
 for the reductive symmetric space ${\mathrm{Sp}}(n,{\mathbb{R}})/{\mathrm{GL}}(n,{\mathbb{R}})$.  
Then using the description of discrete series representations
 \cite{MO84, Vogan88}, 
 we get Counterexample \ref{ex:small} (2).  

We now turn to an analytic approach
 to the question 
 of discrete decomposability
 of the restriction.  
For simplicity, 
 assume $K$ is connected.  
We take a maximal torus $T$ of $K$, 
 and write ${\mathfrak{t}}_0$
 for its Lie algebra.  
Fix a positive system $\Delta^+({\mathfrak{k}}, {\mathfrak{t}})$
and denote by $C_+$ ($\subset \sqrt{-1}{\mathfrak{t}}_0^{\ast}$)
 the dominant Weyl chamber.  
We regard $\widehat T$ 
 as a subset of $\sqrt{-1}{\mathfrak{t}}_0^{\ast}$, 
 and set $\Lambda_+:= C_+ \cap \widehat T$.  
Then Cartan--Weyl highest weight theory
 gives a bijection
\[
   \Lambda_+ \simeq \widehat K, 
   \qquad
   \lambda \mapsto \tau_{\lambda}.  
\]
We recall 
 for a subset $S$ of ${\mathbb{R}}^N$, 
 the asymptotic cone $S\infty$
 is the closed cone
 defined by 
\begin{align*}
  S \infty
:=
  \{ y \in {\mathbb{R}}^N
    :\enspace
&
    \text{there exists a sequence 
    $(y_n, \varepsilon_n) \in S \times {\mathbb{R}}_{>0}$
    such that}
\\
&\text{    $\lim_{n \to \infty} \varepsilon_n y_n = y$
    and 
    $\lim_{n \to \infty} \varepsilon_n =0$}
\}.  
\end{align*}
The asymptotic $K$-support $\operatorname{AS}_K(X)$
 of a $K$-module $X$
 is defined 
 by Kashiwara and Vergne \cite{xkashiv}
 as the asymptotic cone 
 of the highest weights
 of irreducible $K$-modules
 occurring in $X$:
\[
  \operatorname{AS}_K(X) := \operatorname{Supp}_K (X) \infty,   
\]
where $\operatorname{Supp}_K (X)$ is the $K$-support of $X$ given 
 by 
\[
  \operatorname{Supp}_K (X)
  :=
  \{\lambda \in \Lambda_+
   :
   \operatorname{Hom}_K (\tau_{\lambda}, X)\ne \{0\}\}.  
\]

For a closed subgroup $K'$ of $K$, 
 we write ${\mathfrak{k}}_0'$
 for its Lie algebra,
 and regard
 $({\mathfrak{k}}_0')^{\perp}=\operatorname{Ker}(\operatorname{pr}:{\mathfrak{k}}_0^{\ast} \to ({\mathfrak{k}}_0')^{\ast})$
 as a subspace of ${\mathfrak{k}}_0$
 via a $K$-invariant inner product on ${\mathfrak{k}}_0$.  
We set
\[
   C_K(K'):= C_+ \cap \sqrt{-1} \operatorname{Ad}^{\ast}(K)({\mathfrak{k}}_0')^{\perp}.  
\]

An estimate of the singularity spectrum
 of the hyperfunction $K$-character of $X$
 yields a criterion
 of \lq\lq{$K'$-admissibility}\rq\rq\ of $X$
 for a subgroup $K'$ of $K$
 (\cite[Theorem 2.8]{xkAnn98} and \cite{deco-euro}):
\begin{fact}
\label{fact:5.5}
Let $G \supset G'$ be a pair of real reductive linear Lie groups with 
compatible maximal compact subgroups $K \supset K'$, 
 and $X$ an irreducible $({\mathfrak{g}}, K)$-module.  
\begin{enumerate}[{\rm (1)}]
\item  
The following two conditions on the triple $(G,G',X)$
 are equivalent:
\begin{enumerate}[{\rm (i)}]
\item 
$X$ is $K'$-admissible, 
 i.e., 
 $\dim \invHom{K'}{\tau}{X|_{K'}}< \infty$
 for all $\tau \in \dual{K'}{}$.  
\item 
$C_K(K') \cap \operatorname{AS}_K (X) =\{0\}$.  
\end{enumerate}
\item 
If one of (therefore either of) (i) and (ii) is satisfied, 
then $X$ is discretely decomposable 
 as a $({\mathfrak {g}}',K')$-module.  
\end{enumerate}
\end{fact}

\subsection{Classification theory of discretely decomposable pairs}
\label{subsec:discdeco}

We begin with two observations.

First, 
for a Riemannian symmetric pair, that is,
 $(G,G')=(G,K)$ where $G'=K'=K$,
 the restriction $X|_{{\mathfrak {g}}'}$
 is obviously discretely decomposable
 as a $({\mathfrak {g}}',K')$-module
 for any irreducible $({\mathfrak {g}},K)$-module $X$, 
 whereas the reductive pair $(G,G')=({\mathrm{SL}}(n,{\mathbb{C}}),{\mathrm{SL}}(n,{\mathbb{R}}))$
 is an opposite extremal case
 as the restriction $X|_{{\mathfrak {g}}'}$
 is never discretely decomposable
 as a $({\mathfrak {g}}',K')$-module
 for any infinite-dimensional irreducible
 $({\mathfrak {g}},K)$-module
 $X$
 (\cite{xkInvent98}).  
There are also intermediate cases
 such as $(G,G')=({\mathrm{SL}}(n,{\mathbb{R}}),{\mathrm{SO}}(p, n-p))$
 for which the restriction $X|_{{\mathfrak {g}}'}$
 is discretely decomposable
 for some infinite-dimensional irreducible $({\mathfrak {g}},K)$-module
 $X$
 and is not for some other $X$.

Secondly, 
 Harish-Chandra's admissibility theorem 
 \cite{HC} asserts that
\[
  \dim_{\mathbb{C}} \invHom {K} \tau {\pi|_{K}} < \infty
\]
for any $\pi \in \dual G{}$ and $\tau \in \dual K {}$.

This may be regarded
 as a statement 
 for a Riemannian symmetric pair $(G,G')=(G,K)$.  
Unfortunately,
 there is a counterexample 
 to an analogous statement
 for the reductive symmetric pair
 $(G,G')=({\mathrm{SO}}(5,{\mathbb{C}}),{\mathrm{SO}}(3, 2))$, 
 namely,
 we proved in \cite{xkdecoaspm} that 
\[
  \dim_{\mathbb{C}} \invHom {G'} \tau {\pi|_{G'}} = \infty
\quad
\text{for some $\pi \in \dual G {}$ and $\tau \in \dual {G'} {}$}.  
\]
However,
 it is plausible
 \cite[Conjecture A]{xkdecoaspm}
 to have 
 a generalization of Harish-Chandra's admissibility 
 in the category ${\mathcal{H C}}$
 of Harish-Chandra modules
 in the following sense:
\[
     \dim \invHom {{\mathfrak {g}}', K'} {\tau_{K'}} {\pi_K|_{{\mathfrak {g}}'}} < \infty
\]
 for any irreducible $({\mathfrak {g}},K)$-module $\pi_K$
 and irreducible $({\mathfrak {g}}',K')$-module $\tau_{K'}$.

In view of these two observations,
 we consider the following conditions 
 (a) -- (d)
 for a pair of real reductive Lie groups $(G, G')$, 
 and raise a problem:
\begin{problem}
\label{prob:discpair}
Classify the pairs $(G,G')$
 of real reductive Lie groups 
 satisfying the condition (a) below (and also (b), (c) or (d)).  
\end{problem}

\begin{enumerate}[(a)]
\item 
there exist an infinite-dimensional irreducible unitary representation
 $\pi$ of $G$
 and an irreducible unitary representation $\tau$ of $G'$
such that
\[
  0 < \dim \operatorname{Hom}_{{\mathfrak {g}}', K'}(\tau_{K'}, \pi_K|_{{\mathfrak {g}}'})< \infty;
\]
\item 
there exist an infinite-dimensional irreducible unitary representation
 $\pi$ of $G$
 and an irreducible unitary representation $\tau$ of $G'$
such that
\[
  0 < \dim \operatorname{Hom}_{{\mathfrak {g}}', K'}(\tau_{K'}, \pi_K|_{{\mathfrak {g}}'});
\]
\item 
there exist an infinite-dimensional irreducible $({\mathfrak {g}}, K)$-module $X$
 and an irreducible $({\mathfrak {g}}', K')$-module $Y$
such that
\[
  0 < \dim \operatorname{Hom}_{{\mathfrak {g}}', K'}(Y, X|_{{\mathfrak {g}}'})< \infty;
\]
\item 
there exist an infinite-dimensional irreducible $({\mathfrak {g}}, K)$-module $X$
 and an irreducible $({\mathfrak {g}}', K')$-module $Y$
such that
\[
  0 < \dim \operatorname{Hom}_{{\mathfrak {g}}', K'}(Y, X|_{{\mathfrak {g}}'}).  
\]
\end{enumerate}
Obviously we have the following implications:
\begin{alignat*}{3}
&\text{(a)}&&\Rightarrow && \text{(b)}
\\
&\Downarrow && && \Downarrow
\\
&\text{(c)}&&\Rightarrow && \text{(d)}
\end{alignat*}
The vertical (inverse) implications
 (c) $\Rightarrow$ (a)
 and (d) $\Rightarrow$ (b) will mean 
finite-multiplicity results like Harish-Chandra's admissibility theorem.

\vskip 1pc
For symmetric pairs,
 Problem \ref{prob:discpair} has been solved
 in \cite[Theorem 5.2]{xkyosh13}:

\begin{theorem}
\label{thm:discpair}
Let $(G,G')$ be a reductive symmetric pair
 defined by an involutive automorphism $\sigma$
 of a simple Lie group $G$.  
Then the following five conditions 
(a), (b), (c), (d), and 
\begin{equation}
\label{eqn:sigmab}
  \sigma \beta \ne -\beta
\end{equation}
are equivalent.  
Here $\beta$ is the highest noncompact root
 with respect to a \lq\lq{$(-\sigma)$-compatible}\rq\rq\ positive system.  
(See \cite{xkyosh13} for a precise definition.)  
\end{theorem}

\begin{examp}
{\rm{
\begin{enumerate}[{\rm (1)}]
\item 
$\sigma = \theta$ (Cartan involution).  
Then \eqref{eqn:sigmab} is obviously satisfied
 because $\theta \beta = \beta$.  
Needless to say, 
 the conditions (a)--(d) hold when $G'=K$.  
\item 
The reductive symmetric pairs 
$
   (G,G')
$
$
=
$
$({\mathrm{SO}}(p_1+p_2,q), {\mathrm{SO}}(p_1) \times {\mathrm{SO}}(p_2, q))$, 
$({\mathrm{SL}}(2n,{\mathbb{R}}),{\mathrm{Sp}}(n,{\mathbb{C}}))$, 
$({\mathrm{SL}}(2n,{\mathbb{R}}),{\mathbb{T}} \cdot {\mathrm{SL}}(n,{\mathbb{C}}))$
 satisfy \eqref{eqn:sigmab}, 
 and therefore (a)--(d).  
\end{enumerate}
}}
\end{examp}
The classification 
 of irreducible symmetric pairs $(G,G')$
 satisfying one of (therefore all of)
 (a)--(d)
 was given in \cite{xkyosh13}.  
It turns out
 that there are fairly many reductive symmetric pairs
 $(G,G')$
 satisfying the five equivalent conditions
 in Theorem \ref{thm:discpair}
 when $G$ does not carry a complex Lie group structure, 
 whereas there are a few such pairs $(G,G')$
 when $G$ is a complex Lie group.  
As a flavor of the classification, 
 we present a list 
 in this particular case.  
For this, 
 we use the following notation, 
 which is slightly different from that used in the other parts
 of this article.  
Let $G_{\mathbb{C}}$ be a complex simple Lie group, 
 and $G_{\mathbb{R}}$ a real form.  
Take a maximal compact subgroup $K_{\mathbb{R}}$ of $G_{\mathbb{R}}$, 
and let $K_{\mathbb{C}}$ be the complexification
 of $K_{\mathbb{R}}$ in $G_{\mathbb{C}}$.
We denote by ${\mathfrak {g}}$, ${\mathfrak {k}}$, and ${\mathfrak {g}}_{\mathbb{R}}$
 the Lie algebras
 of $G_{\mathbb{C}}$, $K_{\mathbb{C}}$, and $G_{\mathbb{R}}$, 
respectively, 
 and write ${\mathfrak {g}} = {\mathfrak {k}} + {\mathfrak {p}}$
 for the complexified Cartan decomposition.  

\begin{examp}
[{\cite[Corollary 5.9]{xkyosh13}}]
{\rm{
The following five conditions on the pairs $(G_{\mathbb{C}}, G_{\mathbb{R}})$ are equivalent:
\begin{enumerate}[{\rm (i)}]
\item 
$(G_{\mathbb{C}}, K_{\mathbb{C}})$ satisfies (a)
 (or equivalently, (b), (c), or (d)).  
\item 
$(G_{\mathbb{C}}, G_{\mathbb{R}})$ satisfies (a)
 (or equivalently, (b), (c), or (d)).  
\item 
The minimal nilpotent orbit of ${\mathfrak {g}}$
 does not intersect ${\mathfrak {g}}_{\mathbb{R}}$.  
\item 
The minimal nilpotent orbit of ${\mathfrak {g}}$
 does not intersect ${\mathfrak {p}}$.  
\item 
The Lie algebras ${\mathfrak {g}}$, ${\mathfrak {k}}$, 
 and ${\mathfrak {g}}_{\mathbb{R}}$ are given 
 in the following table: 

\medskip
\begin{tabular}{c|ccccc}
  $\mathfrak{g}$ 
& $\mathfrak{sl}(2n,{\mathbb{C}})$ 
& $\mathfrak{so}(m,{\mathbb{C}})$ 
& $\mathfrak{sp}(p+q,{\mathbb{C}})$
&  $\mathfrak{f}_{4}^{{\mathbb{C}}}$ 
& $\mathfrak{e}_{6}^{{\mathbb{C}}}$ 
\\ \hline
$\mathfrak{k}$ 
& $\mathfrak{sp}(n,{\mathbb{C}})$ 
&\,\, $\mathfrak{so}(m-1,{\mathbb{C}})$\,\, 
&\,\, $\mathfrak{sp}(p,{\mathbb{C}}) + \mathfrak{sp}(q,{\mathbb{C}})$\,\,
&\,\,  $\mathfrak{so}(9,{\mathbb{C}})$ \,\,
& $\mathfrak{f}_{4}^{{\mathbb{C}}}$ 
\\ \hline
$\mathfrak{g}_{\mathbb{R}}$ 
& $\mathfrak{su}^{\ast}(2n)$ 
& $\mathfrak{so}(m-1,1)$ 
& $\mathfrak{sp}(p,q)$
&  $\mathfrak{f}_{4(-20)}$ 
&\,\, $\mathfrak{e}_{6(-26)}$ 
\end{tabular}
\medskip
\newline
where $m \ge 5$ and $n,p,q \ge 1$.  
\end{enumerate}
}}
\end{examp}

\begin{rem}
The equivalence (iv) and (v) was announced
 by Brylinski--Kostant
 in the context
 that there is no minimal representation of a Lie group $G_{\mathbb{R}}$
 with the Lie algebra ${\mathfrak{g}}_{\mathbb{R}}$
 in the above table 
(see \cite{Brylinski98}).  
The new ingredient here 
 is that this condition on the Lie algebras
 corresponds to a question of discretely decomposable restrictions
 of Harish-Chandra modules.  
\end{rem}

For nonsymmetric pairs,
 there are a few nontrivial cases
 where (a) (and therefore (b), (c), and (d))
 holds,
 as follows.  

\begin{examp}
[{\cite{xkInvent94}}]
\label{ex:G2}
{\rm{
The nonsymmetric pairs $(G,G')=({\mathrm{SO}}(4,3), {\mathrm{G}}_{2(2)})$
 and 
 $({\mathrm{SO}}(7,{\mathbb{C}}), {\mathrm{G}}_{2}^{{\mathbb{C}}})$ satisfy 
 (a) (and also (b), (c), and (d)).  
}}
\end{examp}

Once we classify the pairs $(G,G')$
 such that there exists at least one irreducible 
 infinite-dimensional $({\mathfrak {g}}, K)$-module $X$
 which is discretely decomposable as a $({\mathfrak {g}}', K')$-module, 
 then we would like to find all such $X$s.

In \cite{decoAq}
 we carried out
 this project for $X=A_{\mathfrak {q}}(\lambda)$
 by applying the general criterion 
 (Facts \ref{nec} and \ref{fact:5.5})
 to reductive symmetric pairs $(G,G')$.  
This is a result in Stage A
 of the branching problem, 
 and we think it will serve as a foundational result
 for Stage B
 (explicit branching laws).  
Here is another example of 
 the classification of the triples $(G,G',X)$
 when $G \simeq G' \times G'$, 
 see {\cite[Theorem 6.1]{xkyosh13}}:

\begin{examp}
[tensor product]
\label{ex:tensor}
{\rm{
Let $G$ be a noncompact connected simple Lie group,
 and let $X_j$ ($j=1,2$) be infinite-dimensional 
irreducible $({\mathfrak {g}}, K)$-modules.  
\begin{enumerate}[{\rm (1)}]
\item 
Suppose $G$ is not of Hermitian type.  
Then the tensor product representation 
 $X_1 \otimes X_2$ is never discretely decomposable
 as a $({\mathfrak {g}}, K)$-module.  
\item 
Suppose $G$ is of Hermitian type.  
Then the tensor product representation $X_1 \otimes X_2$
 is discretely decomposable
 as a $({\mathfrak {g}}, K)$-module
 if and only if 
 both $X_1$ and $X_2$ are simultaneously highest weight
 $({\mathfrak {g}}, K)$-modules
 or simultaneously lowest weight $({\mathfrak {g}}, K)$-modules.  
\end{enumerate}
}}
\end{examp}

\subsection{Two spaces
 $\invHom{G'}{\tau}{\pi|_{G'}}$
 and $\invHom{{\mathfrak {g}}', K'}{\tau_{K'}}{\pi_K|_{\mathfrak{g}'}}$}
There is a canonical injective homomorphism
\[
{\operatorname{Hom}}_{{\mathfrak {g}}', K'}(\tau_{K'}, \pi_K|_{\mathfrak {g}'})
\hookrightarrow
 {\operatorname{Hom}}_{G'}(\tau, \pi|_{G'}), 
\]
however, 
 it is not bijective
 for $\tau \in \widehat{G'}$ and $\pi \in \widehat G$.  
In fact, 
 we have:
\begin{proposition}
\label{prop:gKG}
Suppose that 
$\pi$ is an irreducible unitary representation of $G$.  
 If the restriction $\pi|_{G'}$
 contains a continuous spectrum
 and if an irreducible unitary representation $\tau$
 of $G'$
 appears as an irreducible summand
 of the restriction $\pi|_{G'}$, 
 then we have
\[
{\operatorname{Hom}}_{{\mathfrak {g}}', K'}(\tau_{K'}, \pi_K|_{{\mathfrak {g}}'})=\{0\}
\ne
{\operatorname{Hom}}_{G'}(\tau, \pi|_{G'}).  
\]
\end{proposition}
\begin{qedproof}
If 
$
     \invHom{{\mathfrak {g}}'}{\tau_{K'}}{\pi_K|_{{\mathfrak {g}}'}}
=
     \invHom{{\mathfrak {g}}',K'}{\tau_{K'}}{\pi_K|_{{\mathfrak {g}}'}}
$ were nonzero, 
then the $({\mathfrak {g}}, K)$-module $\pi_K$
 would be discretely decomposable as a $({\mathfrak {g}}', K')$-module
 by Theorem \ref{thm:deco98}.  
In turn, 
 the restriction $\pi|_{G'}$
 of the unitary representation $\pi$ would decompose discretely
 into a Hilbert direct sum 
 of irreducible unitary representations
 of $G'$
 by \cite[Theorem 2.7]{xkdecoaspm}, 
 contradicting the assumption.  
Hence we conclude
$
  {\operatorname{Hom}}_{{\mathfrak {g}}',K'}(\tau_{K'}, \pi_K|_{{\mathfrak {g}}'}) =\{0\}.  
$
\end{qedproof}

An example of Proposition \ref{prop:gKG}
 may be found in \cite[Part II]{xkors}
 where $\pi$ is the minimal representation 
of $G={\mathrm{O}}(p,q)$
 and $\tau$ is the unitarization of a Zuckerman derived functor
 module $A_{\mathfrak {q}}(\lambda)$
 for $G'={\mathrm{O}}(p',q') \times {\mathrm{O}}(p'', q'')$
 with $p=p'+p''$ and $q=q'+q''$
 ($p',q',p'',q'' >1$ and $p+q$ even).

Here is another example of Proposition \ref{prop:gKG}:
\begin{theorem}
\label{thm:tensorP}
Let $G$ be a real reductive linear Lie group,
 and let $\pi= \operatorname{Ind}_P^G ({\mathbb{C}}_{\lambda})$
 be a spherical unitary degenerate principal series representation 
of $G$ 
 induced from a unitary character ${\mathbb{C}}_{\lambda}$
 of a parabolic subgroup $P=L N$ of $G$.  
\begin{enumerate}[{\rm (1)}]
\item 
For any irreducible $({\mathfrak {g}},K)$-module
 $\tau_K$, 
 we have
\[
  {\operatorname{Hom}}_{{\mathfrak {g}},K}
  (\tau_{K}, \pi_K \otimes \pi_K) =\{0\}.  
\]
\item 
Suppose now $G$ is a classical group.  
If $N$ is abelian 
 and $P$ is conjugate 
 to the opposite parabolic subgroup 
 $\overline P = L \overline N$, 
 then we have a unitary equivalence
 of the discrete part:
\begin{equation}
\label{eqn:GLdisc}
  L^2(G/L)_{\operatorname{disc}}
  \simeq
  {\sum_{\tau \in \widehat G}}{}^{\oplus}
  \dim_{\mathbb{C}} \invHom{G}{\tau}{\pi \widehat \otimes \pi} \, \tau.  
\end{equation}
In particular, 
we have
\[
  \dim_{\mathbb{C}} \invHom G{\tau}{\pi \widehat \otimes \pi} \le 1
\]
for any irreducible unitary representation $\tau$ of $G$.  
Moreover there exist countably many irreducible unitary representations $\tau$
 of $G$ such that
\[
  \dim_{\mathbb{C}} \invHom G{\tau}{\pi \widehat \otimes \pi} = 1. 
\]
\end{enumerate}
\end{theorem}
A typical example of the setting
 in Theorem \ref{thm:tensorP} (2) is the Siegel parabolic subgroup
 $P=LN={\mathrm{GL}}(n,{\mathbb{R}}) \ltimes {\mathrm{Sym}}(n,{\mathbb{R}})$
 in $G={\mathrm{Sp}}(n,{\mathbb{R}})$.  

\begin{qedproof}
\begin{enumerate}[(1)]
\item 
This is a direct consequence of Example \ref{ex:tensor}.  
\item 
Take $w_0 \in G$ such that 
 $w_0 L w_0^{-1}  =L$
 and $w_0 N w_0^{-1}  =\overline N$.  
Then the $G$-orbit through $(w_0 P, eP)$
 in $G/P \times G/P$
 under the diagonal action
 is open dense,
 and therefore Mackey theory gives a unitary equivalence
\begin{equation}
\label{eqn:Mackey}
L^2(G/L) \simeq \pi_{\lambda} \widehat \otimes \pi_{\lambda}
\end{equation}
 because $\operatorname{Ad}^{\ast}(w_0) \lambda=-\lambda$, 
see \cite{xkrons}
 for instance.  
Since $N$ is abelian, 
 $(G,L)$ forms a symmetric pair
 (see \cite{xrr92}).  
Therefore the branching law of the tensor product representation 
$\pi \widehat \otimes \pi$
 reduces to the Plancherel formula
 for the regular representation
 on the reductive symmetric space $G/L$, 
 which is known;
 see \cite{xdelorme}.  
In particular,
 we have the unitary equivalence \eqref{eqn:GLdisc}, 
and the left-hand side of \eqref{eqn:GLdisc} is nonzero
 if and only if $\operatorname{rank} G/L =\operatorname{rank} K/L \cap K$
 due to Flensted-Jensen and Matsuki--Oshima \cite{MO84}.  
By the description
 of discrete series representation
 for $G/L$
 by Matsuki--Oshima \cite{MO84}
 and Vogan \cite{Vogan88}, 
 we have the conclusion. 
\end{enumerate}
\end{qedproof}

\subsection{Analytic vectors and discrete decomposability}
Suppose $\pi$ is an irreducible unitary representation of $G$
 on a Hilbert space $V$, 
 and $G'$ is a reductive subgroup of $G$
 as before.  
Any $G'$-invariant closed subspace $W$ in $V$
 contains $G'$-analytic vectors
(hence, 
also $G'$-smooth vectors)
 as a dense subspace.  
However,
 $W$ may not contain nonzero $G$-smooth vectors
 (hence, also $G$-analytic vectors).  
In view of Theorem \ref{thm:deco98}
 in the category ${\mathcal{H C}}$
 of Harish-Chandra modules,
 we think that this is related to the existence
 of a continuous spectrum 
 in the branching law
 of the restriction $\pi|_{G'}$.  
We formulate a problem related to this delicate point below.  
As before,
 $\pi^{\infty}$ and $\tau^{\infty}$ denote
 the space of $G$-smooth vectors
 and $G'$-smooth vectors
 for representations $\pi$ and $\tau$
 of $G$ and $G'$, 
respectively.  
An analogous notation 
 is applied to $\pi^{\omega}$ and $\tau^{\omega}$.  
\begin{problem}
\label{prob:smoothdisc}
Let $(\pi, V)$ be an irreducible unitary representation
 of $G$, 
 and  $G'$ a reductive subgroup of $G$.  
Are the following four conditions
 on the triple $(G,G', \pi)$ equivalent?
\begin{enumerate}[{\rm (i)}]
\item 
There exists an irreducible $({\mathfrak {g}}',K')$-module
 $\tau_{K'}$ such that 
\[
{\operatorname{Hom}}_{{\mathfrak {g}}',K'}(\tau_{K'}, \pi_K|_{\mathfrak{g}'})
\ne 
\{0\}.  
\]
\item 
There exists an irreducible unitary representation $\tau$ of $G'$
 such that 
\[
{\operatorname{Hom}}_{G'}(\tau^{\omega}, \pi^{\omega}|_{G'})
\ne 
\{0\}.  
\]
\item 
There exists an irreducible unitary representation $\tau$ of $G'$
 such that 
\[
{\operatorname{Hom}}_{G'}(\tau^{\infty}, \pi^{\infty}|_{G'})
\ne 
\{0\}.  
\]
\item 
The restriction $\pi|_{G'}$ decomposes
 discretely 
 into a Hilbert direct sum
 of irreducible unitary representations of $G'$.  
\end{enumerate}
\end{problem}
Here are some remarks 
 on Problem \ref{prob:smoothdisc}.  
\begin{rem}
\label{rem:smoothdisc}
{\rm{
\begin{enumerate}[(1)]
\item 
In general,
 the implication (i) $\Rightarrow$ (iv)
 holds (\cite[Theorem 2.7]{xkdecoaspm}).  
\item 
If the restriction $\pi|_{K'}$ is $K'$-admissible,
 then (i) holds by \cite[Proposition 1.6]{xkInvent98}
 and (iv) holds by \cite[Theorem 1.2]{xkInvent94}.   
\item 
The implication (iv) $\Rightarrow$ (i) was raised
 in \cite[Conjecture D]{xkdecoaspm}, 
 and some affirmative results has been announced 
 by Duflo and Vargas 
 in a special setting 
 where $\pi$ is Harish-Chandra's discrete series representation
(cf. \cite{xdv}).  
A related result is given in \cite{xzhuliang}.  
\item  
Even when the unitary representation $\pi|_{G'}$ decomposes discretely 
 (i.e., 
 (iv) in Problem \ref{prob:smoothdisc} holds), 
 it may happen
 that $V^{\infty} \subsetneqq (V|_{G'})^{\infty}$.  
The simplest example
for this is as follows.  
Let $(\pi',V')$ and $(\pi'',V'')$ be infinite-dimensional 
 unitary representations
 of noncompact Lie groups $G'$ and $G''$,
 respectively.  
Set $G=G' \times G''$, 
 with $G'$ realized as a subgroup of $G$
 as $G' \times \{e\}$, 
and set $\pi= \pi' \boxtimes \pi''$.   
Then 
$
   V^{\infty} \subsetneqq
   (V|_{G'})^{\infty}
$
because
$
   (V'')^{\infty}
    \subsetneqq
    V''
$.  
\end{enumerate}
}}
\end{rem}

\section{Features of the restriction, II
: $\invHom{G'}{\pi|_{G'}}{\tau}$ (symmetry breaking operators)}
\label{sec:A2}

In the previous section, 
 we discussed embeddings
 of irreducible $G'$-modules $\tau$
 into irreducible $G$-modules $\pi$
 (or the analogous problem in the category ${\mathcal{HC}}$ 
 of Harish-Chandra modules); 
 see Case I in Section \ref{sec:ep}.  
In contrast, 
 we consider the opposite order in this section, namely,
 continuous $G'$-homomorphisms from irreducible $G$-modules $\pi$
 to irreducible $G'$-modules $\tau$, 
 see Case II in Section \ref{sec:ep}.  
We highlight the case where $\pi$ and $\tau$ are admissible smooth representations
 (Casselman--Wallach globalization
 of modules in the category ${\mathcal{HC}}$).  
Then it turns out
 that the spaces
 $\invHom{G'}{\pi^{\infty}|_{G'}}{\tau^{\infty}}$
 or $\invHom{{\mathfrak {g}}',K'}{\pi_K|_{\mathfrak {g}'}}{\tau_{K'}}$
 are much larger 
 in general than the spaces
 $\invHom{G'}{\tau^{\infty}}{\pi^{\infty}|_{G'}}$
 or $\invHom{{\mathfrak {g}}',K'}{\tau_{K'}}{\pi_K|_{\mathfrak {g}'}}$
 considered in Section \ref{sec:A}.  
Thus the primary concern here will be with obtaining an upper estimate
 for the dimensions of those spaces.

It would make reasonable sense
 to find branching laws (Stage B)
 or to construct symmetry breaking operators (Stage C)
 if we know {\it{a priori}} the nature of the multiplicities
 in branching laws.  
The task of Stage A
 of the branching problem is
 to establish a criterion 
 and to give a classification
 of desirable settings.  
In this section,
 we consider:
\begin{problem}
\label{prob:PPBB}
\begin{enumerate}[{\rm (1)}]
\item 
{\rm{(finite multiplicities)}}\enspace
Find a criterion for when a pair $(G,G')$
 of real reductive Lie groups satisfies
\[
\dim \invHom{G'}{\pi^{\infty}|_{G'}}{\tau^{\infty}}<\infty
\quad
 \text{for any $\pi^{\infty} \in \dual{G}{smooth}
$ and 
               $\tau^{\infty} \in \dual{G'}{smooth}
$.  
}
\]
Classify all such pairs $(G,G')$.  
\item 
{\rm{(uniformly bounded multiplicities)}}\enspace
Find a criterion for when a pair $(G,G')$
 of real reductive Lie groups satisfies
\[
\underset{\pi^{\infty} \in \dual{G}{smooth}}{\sup}
\underset{\tau^{\infty} \in \dual{G'}{smooth}}{\sup}
\dim \invHom{G'}{\pi^{\infty}|_{G'}}{\tau^{\infty}}<\infty.  
\]
Classify all such pairs $(G,G')$.  
\end{enumerate}
\end{problem}
One may also think of variants of Problem \ref{prob:PPBB}.  
For instance,
 we may refine Problem \ref{prob:PPBB}
 by considering it as a condition
 on the triple $(G,G',\pi)$
 instead of a condition on the pair
 $(G,G')$:
\begin{problem}
\label{prob:5.2}
\begin{enumerate}[{\rm (1)}]
\item 
Classify the triples $(G,G',\pi^{\infty})$ 
 with $G \supset G'$
 and $\pi^{\infty} \in \dual {G}{smooth}$ 
such that
\begin{equation}
\label{eqn:3fin}
\dim \invHom{G'}{\pi^{\infty}|_{G'}}{\tau^{\infty}}<\infty
\quad
 \text{for any $\tau^{\infty} \in \dual{G'}{smooth}
$.  
}
\end{equation}
\item 
Classify the triples $(G,G',\pi^{\infty})$ 
 such that
\begin{equation}
\label{eqn:3bdd}
\underset{\tau^{\infty} \in \dual{G'}{smooth}}{\sup}
\dim \invHom{G'}{\pi^{\infty}|_{G'}}{\tau^{\infty}}<\infty.  
\end{equation}
\end{enumerate}
\end{problem}

Problem \ref{prob:PPBB} has been solved recently
 for all reductive symmetric pairs $(G,G')$; 
 see Sections \ref{subsec:BP} and \ref{subsec:classPP}.  
On the other hand, 
 Problem \ref{prob:5.2} has no complete solution
 even when $(G,G')$ is a reductive symmetric pair.  
Here are some partial answers to Problem \ref{prob:5.2} (1):

\begin{examp}
{\rm{
\begin{enumerate}[{\rm (1)}]
\item 
 If $(G,G')$ satisfies (PP)
 (see the list in Theorem \ref{thm:PP}), 
 then the triple $(G,G',\pi)$ satisfies
 \eqref{eqn:3fin} 
 whenever $\pi^{\infty} \in \dual {G}{smooth}$.  
\item 
If $\pi$ is $K'$-admissible, 
 then \eqref{eqn:3fin} is satisfied.  
A necessary and sufficient condition 
 for the $K'$-admissibility
 of $\pi|_{K'}$, 
 Fact \ref{fact:5.5}, 
 is easy to check in many cases.  
In particular, 
 a complete classification 
 of the triples $(G,G',\pi)$
 such that $\pi|_{K'}$
 is $K'$-admissible was recently accomplished in \cite{decoAq}
 in the setting where $\pi_K =A_{\mathfrak {q}}(\lambda)$
 and where $(G,G')$ is a reductive symmetric pair.  
\end{enumerate}
}}
\end{examp}

We give a conjectural statement
 concerning Problem \ref{prob:5.2} (2).  
\begin{conj}
\label{conj:3bdd}
Let $(G,G')$ be a reductive symmetric pair.  
If $\pi$ is an irreducible highest weight representation 
 of $G$
 or if $\pi$ is a minimal representation of $G$, 
 then the uniform boundedness property \eqref{eqn:3bdd}
 would hold 
 for the triple $(G,G',\pi^{\infty})$.  
\end{conj}

Some evidence was given in \cite[Theorems B and D]{K08} 
 and in \cite{xkors, KOP}.  

\subsection{Real spherical homogeneous spaces}
\label{subsec:HP}

A complex manifold $X_{\mathbb{C}}$ 
 with an action of a complex reductive group $G_{\mathbb{C}}$
 is called {\it{spherical}} 
 if a Borel subgroup of $G_{\mathbb{C}}$
 has an open orbit in $X_{\mathbb{C}}$.  
Spherical varieties
 have been studied extensively in the context of algebraic geometry
 and finite-dimensional representation theory.  
In the real setting, 
 in search of a broader framework for global analysis 
 on homogeneous spaces
 than the usual 
 (e.g., 
 reductive symmetric spaces),
 we propose the following:
\begin{definition}
[{\cite{Ksuron}}]
\label{def:realsp}
\rm{
Let $G$ be a real reductive Lie group.  
We say a connected smooth manifold $X$ with $G$-action
 is {\it{real spherical}}
 if a minimal parabolic subgroup $P$
 of $G$ has an open orbit in $X$, 
 or equivalently $\#(P \backslash X)< \infty$.  
}
\end{definition}

The equivalence in Definition \ref{def:realsp}
 was proved in \cite{xbien}
 by using Kimelfeld \cite{Kimelfeld}
 and Matsuki \cite{Matsuki-icm90};
 see \cite[Remark]{xktoshima}
 and references therein for related earlier results.

Here are some partial results
 on the classification 
 of real spherical homogeneous spaces.  
\begin{examp}
\label{ex:HP}
{\rm{
\begin{enumerate}[{\rm (1)}]
\item 
If $G$ is compact 
then all $G$-homogeneous spaces
 are real spherical. 
\item 
Any semisimple symmetric space $G/H$
 is real spherical.  
The (infinitesimal) classification 
 of semisimple symmetric spaces
 was accomplished by Berger \cite{ber}.  
\item 
$G/N$ is real spherical 
 where $N$ is a maximal unipotent subgroup of $G$.  
\item 
For $G$ of real rank one,
 real spherical homogeneous spaces
 of $G$ are classified by Kimelfeld \cite{Kimelfeld}.  
\item 
Any real form $G/H$ 
 of a spherical homogeneous space $G_{\mathbb{C}}/H_{\mathbb{C}}$
 is real spherical \cite[Lemma 4.2]{xktoshima}.  
The latter were classified by 
 Kr{\"a}mer \cite{xkramer}, Brion, \cite{xbrion}, 
 and Mikityuk \cite{xmik}. 
In particular,
 if $G$ is quasi-split,
 then the classification problem 
 of real spherical homogeneous spaces $G/H$ reduces to that
 of the known classification
 of spherical homogeneous spaces.  
\item 
The triple product space $(G \times G \times G)/\operatorname{diag} G$ 
 is real spherical 
 if and only if $G$ is locally isomorphic 
 to the direct product 
 of compact Lie groups
 and some copies of ${\mathrm{O}}(n,1)$
 (Kobayashi \cite{Ksuron}).  
\item 
Real spherical homogeneous spaces
 of the form $(G \times G')/\operatorname{diag} G'$
 for symmetric pairs $(G, G')$ were recently classified.  
We review this in Theorem \ref{thm:PP} below.  
\end{enumerate}
}}
\end{examp}
The second and third examples
 form the basic geometric settings 
 for analysis
 on reductive symmetric spaces
 and Whittaker models.  
The last two examples play a role
 in Stage A 
 of the branching problem, 
 as we see in the next subsection.

The significance of this geometric property
is that the group $G$ controls the space of functions on $X$
 in the sense
 that the finite-multiplicity property
 holds for the regular representation
 of $G$ on $C^{\infty}(X)$:
\begin{fact}
[{\cite[Theorems A and C]{xktoshima}}]
\label{fact:HP}
Suppose $G$ is a real reductive linear Lie group,
 and $H$ is an algebraic reductive subgroup. 
\begin{enumerate}[{\rm (1)}]
\item  
The homogeneous space $G/H$ 
 is real spherical
 if and only if 
\[
\text{
$\invHom{G}{\pi^{\infty}}{C^{\infty}(G/H)}$
 is finite-dimensional
 for all $\pi^{\infty} \in \dual G {smooth}$.  
}
\]
\item 
The complexification $G_{\mathbb{C}}/H_{\mathbb{C}}$
 is spherical 
 if and only if 
\[
\sup_{\pi^{\infty} \in \dual G {smooth}}
 \dim_{\mathbb{C}} \invHom G{\pi^{\infty}}{C^{\infty}(G/H)}
< \infty.  
\]
\end{enumerate}
\end{fact}
See \cite{xktoshima}
 for upper and lower estimates
 of the dimension,
 and also for the non-reductive case.  
The proof uses the theory 
 of regular singularities
 of a system of partial differential equations
 by taking an appropriate compactification
 with normal crossing boundaries.  

\subsection{A geometric estimate
 of multiplicities : (PP) and (BB)}
\label{subsec:PP}

Suppose that $G'$ is an algebraic reductive subgroup of 
 $G$.  
For Stage A in the branching problem
 for the restriction $G \downarrow G'$, 
 we apply the general theory
 of Section \ref{subsec:HP}
 to the homogeneous space
 $(G \times G')/\operatorname {diag} G'$.

Let $P$ be a minimal parabolic subgroup of $G$, 
 and $P'$ a minimal parabolic subgroup
 of $G'$.  

\begin{def-lemma}
[{\cite{xktoshima}}]
\label{def:pp}
We say the pair $(G,G')$ satisfies the property (PP)
 if one of the following five equivalent
 conditions is satisfied:  
\par\noindent
{\rm{(PP1)}}\enspace
$(G \times G')/\operatorname{diag}G'$ is real spherical
 as a $(G \times G')$-space.  
\par\noindent
{\rm{(PP2)}}\enspace
$G/P'$ is real spherical as a $G$-space.  
\par\noindent
{\rm{(PP3)}}\enspace
$G/P$ is real spherical as
 a $G'$-space. 
\par\noindent
{\rm{(PP4)}}\enspace
$G$ has an open orbit in $G/P \times G/P'$
 via the diagonal action.  
\par\noindent
{\rm{(PP5)}}\enspace
$\# (P' \backslash G/P) < \infty$.  
\end{def-lemma}
Since the above five equivalent conditions are determined 
 by the Lie algebras ${\mathfrak {g}}$ and ${\mathfrak {g}}'$, 
 we also say
 that the pair $({\mathfrak {g}}, {\mathfrak {g}}')$
 of reductive Lie algebras
 satisfies (PP), 
 where ${\mathfrak {g}}$ and ${\mathfrak {g}}'$ are the Lie algebras
 of the Lie groups $G$ and $G'$, 
 respectively.  
\begin{rem}
\label{rem:PP5}
{\rm{
If the pair $({\mathfrak {g}}, {\mathfrak {g}}')$ satisfies (PP), 
 in particular, 
 (PP5), 
 then there are only finitely many possibilities
 for $\operatorname{Supp} T$
 for symmetry breaking operators
 $T:C^{\infty}(G/P, {\mathcal {V}})$ $\to C^{\infty}(G'/P', {\mathcal {W}})$
 (see Definition \ref{def:SuppT} below).  
This observation has become a guiding principle
 to formalise a strategy
 in classifying all symmetry breaking operators used in \cite{xtsbon}, 
 as we shall discuss in Section \ref{subsec:dist}  
}}
\end{rem}

\vskip 1pc
Next we consider another property, 
 to be denoted (BB), 
 which is stronger than (PP).  
Let $G_{\mathbb{C}}$ be a complex Lie group
 with Lie algebra 
 ${\mathfrak{g}}_{\mathbb{C}}={\mathfrak{g}}\otimes_{\mathbb{R}}{\mathbb{C}}$, 
 and $G_{\mathbb{C}}'$ a subgroup of $G_{\mathbb{C}}$
 with complexified Lie algebra 
 ${\mathfrak{g}}_{\mathbb{C}}'
 ={\mathfrak{g}}'\otimes_{\mathbb{R}}{\mathbb{C}}$.  
We do not assume
 either $G \subset G_{\mathbb{C}}$ or $G' \subset G_{\mathbb{C}}'$.  
Let $B_{\mathbb{C}}$ and $B_{\mathbb{C}}'$ be Borel subgroups of $G_{\mathbb{C}}$ 
 and $G_{\mathbb{C}}'$, 
 respectively.  

\begin{def-lemma}
\label{def:BB}
We say the pair
 $(G,G')$
 (or the pair $({\mathfrak{g}}, {\mathfrak{g}}')$)
 satisfies the property {\rm{(BB)}}
 if one of the following five equivalent conditions is satisfied:
\par\noindent
{\rm{(BB1)}}\enspace
$(G_{\mathbb{C}} \times G_{\mathbb{C}}')/
\operatorname{diag}G_{\mathbb{C}}'$
 is spherical 
 as a $(G_{\mathbb{C}} \times G_{\mathbb{C}}')$-space.   
\par\noindent
{\rm{(BB2)}}\enspace
$G_{\mathbb{C}}/B_{\mathbb{C}}'$ is spherical 
 as a $G_{\mathbb{C}}$-space.  
\par\noindent
{\rm{(BB3)}}\enspace
 $G_{\mathbb{C}}/B_{\mathbb{C}}$ is spherical as
 a $G_{\mathbb{C}}'$-space.
\par\noindent
{\rm{(BB4)}}\enspace
$G_{\mathbb{C}}$ has an open orbit in $G_{\mathbb{C}}/B_{\mathbb{C}} \times G_{\mathbb{C}}/B_{\mathbb{C}}'$
 via the diagonal action.  
\par\noindent
{\rm{(BB5)}}\enspace
$\#(B_{\mathbb{C}}' \backslash G_{\mathbb{C}}/B_{\mathbb{C}}) < \infty$.  
\end{def-lemma}
The above five equivalent conditions (BB1) -- (BB5)
 are determined 
 only by the complexified Lie algebras
 ${\mathfrak {g}}_{\mathbb{C}}$
 and 
 ${\mathfrak {g}}_{\mathbb{C}}'$.  

\begin{rem}
\label{rem:PPBB}
\begin{enumerate}[{\rm (1)}]
\item 
{\rm{(BB)}} implies {\rm{(PP)}}.  
\item 
If both $G$ and $G'$ are quasi-split,
 then {\rm{(BB)}} $\Leftrightarrow$ {\rm{(PP)}}.  
\end{enumerate}
In fact, 
the first statement follows immediately from \cite[Lemmas 4.2 and 5.3]{xktoshima}, 
 and the second statement is clear.  

\end{rem}

\subsection{Criteria for finiteness/boundedness
 of multiplicities}
\label{subsec:BP}
In this and the next subsections, 
 we give an answer to Problem \ref{prob:PPBB}.  
The following criteria are direct consequences
 of Fact \ref{fact:HP}
 and a careful consideration of the topology 
 of representation spaces, 
 and are proved in \cite{xktoshima}.  
\begin{theorem}
\label{thm:A}
The following three conditions
on a pair of real reductive
 algebraic groups $G\supset G'$
 are equivalent:
\begin{enumerate}[{\rm (i)}]
\item 
{\rm{(Symmetry breaking)}}\enspace
$\invHom{G'}{\pi^{\infty}|_{G'}}{\tau^{\infty}}$
 is finite-dimensional 
 for any pair $({\pi^{\infty}},{\tau^{\infty}})$
 of irreducible smooth representations
 of $G$ and $G'$.  
\item 
{\rm{(Invariant bilinear form)}}\enspace
There exist at most finitely many 
 linearly independent
 $G'$-invariant bilinear forms
 on $\pi^{\infty}|_{G'} \widehat \otimes \tau^{\infty}$, 
 for any $\pi^{\infty} \in \dual {G}{smooth}$
 and $\tau^{\infty} \in \dual {G'}{smooth}$.  
\item 
{\rm{(Geometry)}}\enspace
The pair $(G, G')$ satisfies the condition (PP)
 (Definition-Lemma \ref{def:pp}).  
\end{enumerate}
\end{theorem}
\begin{theorem}
\label{thm:B}
The following three conditions
 on a pair of real reductive algebraic groups
 $G \supset G'$
 are equivalent:
\begin{enumerate}[{\rm (i)}]
\item
{\rm{(Symmetry breaking)}}\enspace
There exists a constant $C$
 such that
\[
\dim_{\mathbb{C}}
\invHom{G'}{\pi^{\infty}|_{G'}}{\tau^{\infty}}
\le C
\]
for any $\pi^{\infty} \in \dual {G}{smooth}$
 and $\tau^{\infty} \in \dual {G'}{smooth}$.  
\item
{\rm{(Invariant bilinear form)}}\enspace
There exists a constant $C$
 such that 
\[
\dim_{\mathbb{C}}
\invHom{G'}{\pi^{\infty}|_{G'} \widehat \otimes \tau^{\infty}}
{{\mathbb{C}}}\le C
\]
for any $\pi^{\infty} \in \dual {G}{smooth}$
 and $\tau^{\infty} \in \dual {G'}{smooth}$.  
\item
{\rm{(Geometry)}}\enspace
The pair $(G,G')$ satisfies the condition (BB)
 (Definition-Lemma \ref{def:BB}).  
\end{enumerate}
\end{theorem}

\subsection{Classification theory of finite-multiplicity branching laws}
\label{subsec:classPP}

This section gives a complete list 
 of the reductive symmetric pairs
 $(G,G')$
 such that 
 $\dim \invHom{G'}{\pi^{\infty}|_{G'}}{\tau^{\infty}}$
 is finite or bounded
 for all $\pi^{\infty} \in \widehat {G}_{\operatorname{smooth}}$
 and $\tau^{\infty} \in \widehat {G'}_{\operatorname{smooth}}$.  
Owing to the criteria in 
 Theorems \ref{thm:A} and \ref{thm:B}, 
 the classification is reduced to that 
 of (real) spherical homogeneous spaces
 of the form $(G \times G')/\operatorname{diag}G'$, 
 which was accomplished in \cite{xKMt}
 by using an idea of \lq\lq{linearization}\rq\rq\ :

\begin{theorem}
\label{thm:PP}
Suppose $(G,G')$ is a reductive symmetric pair.  
Then the following two conditions
 are equivalent:
\begin{enumerate}[{\rm (i)}]
\item 
$\invHom{G'}{\pi^{\infty}|_{G'}}{\tau^{\infty}}$ 
 is finite-dimensional 
 for any pair $(\pi^{\infty}, \tau^{\infty})$
 of admissible smooth representations of $G$ and $G'$.  
\item 
The pair $({\mathfrak{g}},{\mathfrak{g}}')$ of their Lie algebras 
 is isomorphic 
 {\rm{(}}up to outer automorphisms{\rm{)}}
 to a direct sum of the following pairs:
\begin{enumerate}
\item[{\rm{A)}}]
{\rm{Trivial case:}}
${\mathfrak{g}}={\mathfrak{g}}'$.  
\item[{\rm{B)}}]
{\rm{Abelian case:}}
${\mathfrak {g}}={\mathbb{R}}$, 
${\mathfrak {g}}'=\{0\}$.  
\item[{\rm{C)}}]
{\rm{Compact case:}}
${\mathfrak {g}}$ is the Lie algebra
 of a compact simple Lie group.  
\item[{\rm{D)}}]
{\rm{Riemannian symmetric pair:}}
${\mathfrak {g}}'$ is the Lie algebra 
 of a maximal compact subgroup $K$
 of a noncompact simple Lie group $G$.  
\item[{\rm{E)}}]
{\rm{Split rank one case ($\operatorname{rank}_{{\mathbb{R}}}G=1$):}}
\begin{enumerate}
\item [{\rm{E1)}}]
$({\mathfrak{o}}(p+q,1),
{\mathfrak{o}}(p)+{\mathfrak{o}}(q,1))$
\quad\quad\,\, $(p+q \ge 2)$, 
\item [{\rm{E2)}}]
$({\mathfrak{su}}(p+q,1),
{\mathfrak{s}}({\mathfrak {u}}(p)+{\mathfrak{u}}(q,1)))$
\quad
$(p+q \ge 1)$, 
\item [{\rm{E3)}}]
$({\mathfrak{sp}}(p+q,1),
{\mathfrak{sp}}(p)+{\mathfrak{sp}}(q,1))$
\quad\,
$(p+q \ge 1)$, 
\item [{\rm{E4)}}]
$({\mathfrak{f}}_{4(-20)},
{\mathfrak{o}}(8,1))$.  
\end{enumerate}
\item[{\rm{F)}}]
{\rm{Strong Gelfand pairs
 and their real forms:}}
\begin{enumerate}
\item[{\rm{F1)}}]
$({\mathfrak{sl}}(n+1,{\mathbb{C}}),
{\mathfrak{gl}}(n,{\mathbb{C}}))$
\, $(n\ge 2)$,  
\item[{\rm{F2)}}]
$({\mathfrak{o}}(n+1,{\mathbb{C}}),
{\mathfrak{o}}(n,{\mathbb{C}}))$
\,\,\,\, $(n\ge 2)$, 
\item[{\rm{F3)}}]
$({\mathfrak{sl}}(n+1,{\mathbb{R}}),
{\mathfrak{gl}}(n,{\mathbb{R}}))$ 
\,\,\,$(n\ge 1)$, 
\item[{\rm{F4)}}]
$({\mathfrak{su}}(p+1,q),{\mathfrak{u}}(p,q))$
\,\,\,\,\, $(p+q\ge 1)$,   
\item[{\rm{F5)}}]
$({\mathfrak{o}}(p+1,q),{\mathfrak{o}}(p,q))$ 
\,\,\,\,\,\,\,\,\,\,$(p+q\ge 2)$.  
\end{enumerate}
\item[{\rm{G)}}]
{\rm{Group case:}}
$({\mathfrak{g}}, {\mathfrak{g}}')=
 ({\mathfrak{g}}_1+{\mathfrak{g}}_1, \operatorname{diag} {\mathfrak{g}}_1)$
where 
\begin{enumerate}
\item[{\rm{G1)}}]
${\mathfrak{g}}_1$ is the Lie algebra
 of a compact simple Lie group, 
\item[{\rm{G2)}}]
$({\mathfrak{o}}(n,1)+{\mathfrak{o}}(n,1), \operatorname{diag}
{\mathfrak{o}}(n,1))$
\quad
$(n \ge 2)$.  
\end{enumerate}
\item[{\rm{H)}}]
{\rm{Other cases:}}
\begin{enumerate}
\item[{\rm{H1)}}]
$({\mathfrak{o}}(2n, 2),
{\mathfrak{u}}(n,1))$
\hphantom{MMMMMMMMM}
$(n \ge 1)$.
\item[{\rm{H2)}}]
$({\mathfrak{su}}^{\ast}(2n+2),
{\mathfrak{su}}(2)+{\mathfrak{su}}^{\ast}(2n)+\mathbb{R})$ 
\quad
$(n\ge 1)$.  
\item[{\rm{H3)}}]
$({\mathfrak{o}}^{\ast}(2n+2),
{\mathfrak{o}}(2)+{\mathfrak{o}}^{\ast}(2n))$
\hphantom{MMMMn} $(n\ge 1)$.
\item[{\rm{H4)}}]
$({\mathfrak{sp}}(p+1,q),
{\mathfrak{sp}}(p,q)+{\mathfrak{sp}}(1))$.  
\item[{\rm{H5)}}]
$({\mathfrak{e}}_{6(-26)},
{\mathfrak{so}}(9,1)+{\mathbb{R}})$.
\end{enumerate}
\end{enumerate}

\end{enumerate}
\end{theorem}

Among the pairs $({\mathfrak {g}}, {\mathfrak {g}}')$
 in the list (A)--(H)
 in Theorem \ref{thm:PP}
 describing finite multiplicities, 
 those pairs
 having uniform bounded multiplicities are classified as follows.  
\begin{theorem}
\label{thm:BB}
Suppose $(G,G')$ is a reductive symmetric pair.  
Then the following two conditions are equivalent:

\begin{enumerate}[{\rm (i)}]
\item 
There exists a constant $C$
 such that 
\[
  \dim_{\mathbb{C}}
 \invHom{G'}{\pi^{\infty}|_{G'}}{\tau^{\infty}}
 \le C
\]
 for any $\pi^{\infty} \in \dual {G}{smooth}$
 and $\tau^{\infty} \in \dual {G'}{smooth}$.
\item 
The pair of their Lie algebras 
 $({\mathfrak{g}},{\mathfrak{g}}')$
 is isomorphic 
 {\rm{(}}up to outer automorphisms{\rm{)}}
 to a direct sum of the pairs in {\rm{(A)}}, {\rm{(B)}}
 and {\rm{(F1)}} -- {\rm{(F5)}}.  
\end{enumerate}
\end{theorem}

\begin{qedproof}
Theorem \ref{thm:PP} follows directly from Theorem \ref{thm:A}
 and \cite[Theorem 1.3]{xKMt}.  
Theorem \ref{thm:BB} follows directly from Theorem \ref{thm:B}
 and \cite[Proposition 1.6]{xKMt}.  
\end{qedproof}

\begin{examp}
\label{ex:AA}
{\rm{
In connection with branching problems,
 some of the pairs appeared earlier
 in the literature.  
For instance,
\begin{alignat*}{3}
&\text{
(F1), (F2)
}
&&
\cdots
&&
\text{
finite-dimensional representations
(strong Gelfand pairs)
\cite{xkramer};
}
\\
&\text{
(F2), (F5)
}
&&
\cdots
&&
\text{
tempered unitary representations
(Gross--Prasad conjecture) 
\cite{GP};
}
\\
&\text{
(G2)
}
&&
\cdots
&&
\text{
tensor product, 
trilinear forms
\cite{CKOP, Ksuron};
}
\\
&\text{
(F1)--(F5)}
&&
\cdots
&&
\text{
multiplicity-free restrictions
\cite{aizenbudgourevitch, xsunzhu}.    
}
\end{alignat*}
}}
\end{examp}

\section{Construction of symmetry breaking operators}
\label{sec:C}
Stage C in the branching problem asks for an explicit construction
of intertwining operators.  
This problem depends on the geometric models
 of representations of a group $G$ and 
 its subgroup $G'$.  
In this section 
 we discuss symmetry breaking operators
 in two models,
 i.e., 
 in the setting of real flag manifolds
 (Sections \ref{subsec:diff}--\ref{subsec:finDiff})
 and in the holomorphic setting
 (Sections \ref{subsec:local}--\ref{subsec:loc}).  

\subsection{Differential operators on different base spaces}
\label{subsec:diff}
We extend the usual notion
 of differential operators
 between two vector bundles
 on the {\it{same}} base space
 to those on {\it{different}} base spaces $X$ and $Y$
 with a morphism $p:Y \to X$ as follows.  
\begin{definition}
\label{def:Peetre}
{\rm{
Let ${\mathcal{V}} \to X$
 and ${\mathcal{W}} \to Y$
 be two vector bundles,
 and $p:Y \to X$
 a smooth map between the base manifolds.  
A continuous linear map 
 $T:C^{\infty}(X,{\mathcal{V}}) \to C^{\infty}(Y,{\mathcal{W}})$
 is said to be a {\it{differential operator}}
 if 
\begin{equation}
\label{eqn:local}
p({\operatorname{Supp}}(T f)) \subset {\operatorname{Supp}} f
\quad
\text{for all } f \in C^{\infty}(X,{\mathcal{V}}), 
\end{equation}
where ${\operatorname{Supp}}$ stands for the support
 of a section.  
}}
\end{definition}
The condition \eqref{eqn:local} shows
that $T$ is a local operator
 in the sense 
 that for any open subset $U$ of $X$, 
 the restriction $(Tf)|_{p^{-1}(U)}$ is determined
 by the restriction $f|_U$.

\begin{examp}
\label{ex:local}
{\rm{
\begin{enumerate}[{\rm (1)}]
\item 
If $X=Y$ and $p$ is the identity map, 
 then the condition \eqref{eqn:local} is equivalent
 to the condition that $T$ is a differential operator
 in the usual sense, 
 due to Peetre's theorem \cite{Pee59}.  
\item 
If $p:Y \to X$ is an immersion,
 then any operator $T$ satisfying \eqref{eqn:local}
 is locally of the form
\[
  \left.
  \sum_{(\alpha, \beta)\in {\mathbb{N}}^{m+n}}
  g_{\alpha \beta}(y) \frac{\partial^{|\alpha|+|\beta|}}{\partial y^{\alpha} \partial z^{\beta}}
\right|_{z_1 = \cdots = z_n =0}
\quad
\text{(finite sum)}, 
\]
where $\{(y_1, \dots, y_m, z_1, \dots,z_n)\}$ are local coordinates
 of $X$
 such that $Y$ is given locally by the equation
 $z_1 = \cdots= z_n=0$, 
and $g_{\alpha \beta}(y)$ are matrix-valued functions on $Y$.  
\end{enumerate}
}}
\end{examp}

\subsection{Distribution kernels
 for symmetry breaking operators}
\label{subsec:dist}
In this section,
 we discuss symmetry breaking operators
 in a geometric setting, 
 where representations are realized
 in the space of smooth sections
 for homogeneous vector bundles.

Let $G$ be a Lie group, 
 and ${\mathcal{V}} \to X$ a homogeneous vector bundle, 
namely,
 a $G$-equivariant vector bundle
 such that the $G$-action on the base manifold $X$ 
 is transitive.  
Likewise,
 let ${\mathcal{W}} \to Y$ be a homogeneous vector bundle
 for a subgroup $G'$.  
The main assumption
 of our setting is 
that there is a $G'$-equivariant map $p:Y \to X$.  
For simplicity,
 we also assume that $p$ is injective, 
 and do not assume any relationship 
 between $p^{\ast} {\mathcal{V}}$ and ${\mathcal{W}}$.  
Then we have continuous representations 
 of $G$ 
 on the Fr{\'e}chet space 
 $C^{\infty}(X, {\mathcal{V}})$
 and of the subgroup $G'$ on $C^{\infty}(Y, {\mathcal{W}})$, 
but it is not obvious
 if there exists a nonzero continuous $G'$-homomorphism
 (symmetry breaking operator)
\[
T: C^{\infty}(X, {\mathcal{V}}) \to C^{\infty}(Y, {\mathcal{W}}).  
\]

In this setting, 
 a basic problem is:
\begin{problem}
\label{prob:tabc}
\begin{enumerate}[{\rm (1)}]
\item 
{\rm{(Stage A)}}\enspace
Find an upper and lower estimate
 of the dimension of the space 

$\invHom{G'}{C^{\infty}(X, {\mathcal{V}})}
{C^{\infty}(Y, {\mathcal{W}})}$
 of symmetry breaking operators.  

\item 
{\rm{(Stage A)}}\enspace
When is $\invHom{G'}{C^{\infty}(X, {\mathcal{V}})}
{C^{\infty}(Y, {\mathcal{W}})}$ finite-dimensional 
 for any $G$-equivariant vector bundle ${\mathcal{V}} \to X$
 and any $G'$-equivariant vector bundle ${\mathcal{W}} \to Y$?  

\item 
{\rm{(Stage B)}}\enspace
Given equivariant vector bundles
 ${\mathcal{V}} \to X$ 
 and ${\mathcal{W}} \to Y$, 
 determine the dimension
 of $\invHom{G'}{C^{\infty}(X, {\mathcal{V}})}
{C^{\infty}(Y, {\mathcal{W}})}$.  

\item 
{\rm{(Stage C)}}\enspace
Construct explicit elements
 in $\invHom{G'}{C^{\infty}(X, {\mathcal{V}})}
{C^{\infty}(Y, {\mathcal{W}})}$.  
\end{enumerate}
\end{problem}

Here are some special cases:
\begin{examp}
\label{ex:GG}
{\rm{
Suppose $G=G'$, 
 $X$ is a (full) real flag manifold $G/P$
 where $P$ is a minimal parabolic subgroup of $G$, 
 and $Y$ is algebraic.  
\begin{enumerate}[{\rm (1)}]
\item 
In this setting, 
 Problem \ref{prob:tabc} (1) and (2) were solved in \cite{xktoshima}.  
In particular, 
 a necessary and sufficient condition
 for Problem \ref{prob:tabc} (2) is 
that $Y$ is real spherical, 
 by Fact \ref{fact:HP} (1)
 (or directly from the original proof of \cite[Theorem A]{xktoshima}).  

\item 
Not much is known 
about precise results
 for Problem \ref{prob:tabc} (3), 
 even when $G=G'$.  
On the other hand, 
 Knapp--Stein intertwining operators
 or Poisson transforms
 are examples
 of explicit intertwining operators
 when $Y$ is a real flag manifold
 or a symmetric space, 
respectively, 
giving a partial solution to Problem \ref{prob:tabc} (4).  
\end{enumerate}
}}
\end{examp}

\begin{examp}
{\rm{
Let $G$ be the conformal group of the standard sphere $X=S^n$, 
 let $G'$ be the subgroup 
 that leaves the totally geodesic submanifold $Y=S^{n-1}$ invariant, 
 and let ${\mathcal{V}} \to X$, ${\mathcal{W}} \to Y$ be 
 $G$-, $G'$-equivariant line bundles, 
respectively.  
Then ${\mathcal{V}}$ and ${\mathcal{W}}$ are parametrized
 by complex numbers $\lambda$ and $\nu$, 
 respectively,
 up to signatures.  
In this setting Problem \ref{prob:tabc} (3) and (4)
 were solved in \cite{xtsbon}.  
This is essentially the geometric setup
 for the classification 
 of $\invHom{{\mathrm{O}}(n,1)}{I(\lambda)^{\infty}}{J(\nu)^{\infty}}$
 which was discussed
 in Section \ref{subsec:2.2}.  
}}
\end{examp}

\vskip 0.5pc
We return to the general setting.  
Let $H$ be an algebraic subgroup of $G$, 
 $(\lambda, V)$ a finite-dimensional representation of $H$, 
 and ${\mathcal {V}}:=G \times_H V \to X:=G/H$
 the associated $G$-homogeneous bundle.  
Likewise, 
let $(\nu, W)$ be a finite-dimensional representation of $H':=H \cap G'$, 
 and ${\mathcal {W}}:=G' \times_{H'} W \to Y:=G'/H'$
 the associated $G'$-equivariant bundle.   
Denote by ${\mathbb{C}}_{2 \rho}$ the one-dimensional representation
 of $H$ defined by 
$
 h \mapsto |\det(\operatorname{Ad}(h):
{\mathfrak {g}}/{\mathfrak {h}}
\to 
{\mathfrak {g}}/{\mathfrak {h}}
)|^{-1}. 
$
Then the volume density bundle $\Omega_{G/H}$ of $G/H$ 
 is given as a homogeneous bundle $G \times_H {\mathbb{C}}_{2 \rho}$.  
Let $(\lambda^{\vee}, V^{\vee})$ be the contragredient representation
 of the finite-dimensional representation $(\lambda, V)$ of $H$.  
Then the dualizing bundle 
 ${\mathcal{V}}^{\ast} := {\mathcal{V}}^{\vee} \otimes \Omega_{G/H}$
is given by 
$
{\mathcal{V}}^{\ast}\simeq G\times_H
 (V^{\vee}
 \otimes {\mathbb{C}}_{2\rho})
$
 as a homogeneous vector bundle.  

\vskip 1pc
By the Schwartz kernel theorem, 
 any continuous operator
 $T:C^{\infty}(X,{\mathcal{V}}) \to C^{\infty}(Y,{\mathcal{W}})$
 is given by a distribution kernel 
$k_T \in {\mathcal{D}}'(X \times Y,{\mathcal{V}}^{\ast} \boxtimes 
{\mathcal{W}})$.  
We write 
\[
 m: G \times G' \to G, 
\quad
 (g,g') \mapsto (g')^{-1}g,  
\]
for the multiplication map.  
If $T$ intertwines $G'$-actions, 
then $k_T$ is $G'$-invariant 
 under the diagonal action, 
 and therefore $k_T$ is of the form $m^{\ast}K_T$
 for some $K_T \in {\mathcal{D}}'(X, {\mathcal {V}}^{\ast}) \otimes W$.  
We have shown in \cite[Proposition 3.1]{xtsbon} 
 the following proposition:  
\begin{proposition}
\label{prop:KT}
Suppose $X$ is compact.  
Then the correspondence
 $T \mapsto K_T$ induces a bijection:
\[
  \operatorname{Hom}_{G'}(C^{\infty}(X,{\mathcal{V}}),C^{\infty}(Y,{\mathcal{W}}))
\overset{\sim}\to
({\mathcal{D}}'(X,{\mathcal{V}}^{\ast})\otimes W)^{\Delta(H')}.  
\]
\end{proposition}
Using Proposition \ref{prop:KT}, 
 we can give a solution to Problem \ref{prob:tabc} (2)
 when $X$ is a real flag manifold:
\begin{theorem}
\label{thm:fmXY}
Suppose $P$ is a minimal parabolic subgroup of $G$, 
 $X=G/P$, 
 and $Y=G'/(G' \cap P)$.  
Then $\operatorname{Hom}_{G'}(C^{\infty}(X,{\mathcal{V}}),C^{\infty}(Y,{\mathcal{W}}))$
is finite-dimensional
 for any $G$-equivariant vector bundle
 ${\mathcal{V}} \to X$
 and any $G'$-equivariant vector bundle
 ${\mathcal{W}} \to Y$
 if and only if 
 $G/(G' \cap P)$
 is real spherical.  
\end{theorem}
\begin{qedproof}
We set $\widetilde Y:=G/(G' \cap P)$
 and $\widetilde {\mathcal{W}}:= G \times_{(G' \cap P)}W$.  
Then Proposition \ref{prop:KT}
 implies that there is a canonical bijection:
\[
\operatorname{Hom}_{G}(C^{\infty}(X,{\mathcal{V}}),C^{\infty}(\widetilde Y,\widetilde {\mathcal{W}}))\overset \sim \to 
\operatorname{Hom}_{G'}(C^{\infty}(X,{\mathcal{V}}),C^{\infty}(Y,{\mathcal{W}})
).  
\]
We apply \cite[Theorem A]{xktoshima}
 to the left-hand side, 
 and get the desired conclusion 
for the right-hand side.  
\end{qedproof}
The smaller $X$ is, 
 the more likely it will be 
 that there exists $Y$ 
 satisfying the finiteness condition
 posed in Problem \ref{prob:tabc} (2).  
Thus one might be interested 
 in replacing the {\it{full}} real flag manifold
 by a {\it{partial}} real flag manifold
 in Theorem \ref{thm:fmXY}.  
By applying the same argument
 as above
 to a generalization of \cite{xktoshima}
 to a partial flag manifold
 in \cite[Corollary 6.8]{xkProg2014}, 
 we get 
\begin{proposition}
\label{prop:fmP}
Suppose $P$ is a (not necessarily minimal) parabolic subgroup 
 of $G$
 and $X=G/P$.  
Then the finiteness condition
 for symmetry breaking operators
 in Problem \ref{prob:tabc} (2)
 holds only 
 if the subgroup $G' \cap P$ has an open orbit in $G/P$.  
\end{proposition}
Back to the general setting, 
 we endow the double coset space
 $H' \backslash G/H$ 
 with the quotient topology 
 via the canonical quotient 
 $G \to H' \backslash G/H$.  
Owing to Proposition \ref{prop:KT}, 
 we associate a closed subset of 
 $H' \backslash G/H$
 to each symmetry breaking operator: 
\begin{definition}
\label{def:SuppT}
{\rm{
Given a continuous symmetry breaking operator
 $T:C^{\infty}(X,{\mathcal{V}}) \to C^{\infty}(Y,{\mathcal{W}})$, 
 we define a closed subset $\operatorname{Supp}T$
 in the double coset space $H' \backslash G/H$
 as the support
 of $K_T \in {\mathcal{D}}'(X,{\mathcal{V}}^{\ast}) \otimes W$.  
}}
\end{definition}

\begin{examp}
{\rm{
If $H=P$, 
 a minimal parabolic subgroup of $G$, 
 and if $H'$ has an open orbit in $G/P$, 
 then $\#(H' \backslash G/P)< \infty$.  
In particular,
 there are only finitely many possibilities
 for ${\operatorname{Supp}}\ T$.  
}}
\end{examp}

\begin{definition}
\label{defprop:regdiff}
{\rm{
Let $T: C^{\infty}(X, {\mathcal{V}}) \to C^{\infty}(Y, {\mathcal{W}})$ be a continuous symmetry breaking operator.  
\begin{enumerate}[1)]
\item 
We say $T$ is a {\it{regular}} symmetry breaking operator
 if $\operatorname{Supp}T$ contains an interior point
 of $H' \backslash G/H$.  
We say $T$ is {\it{singular}}
 if $T$ is not regular.  
\item 
We say $T$ is a {\it{differential}} symmetry breaking operator
if $\operatorname{Supp}T$ is a singleton in $H' \backslash G/H$.  
\end{enumerate}
}}
\end{definition}

\begin{rem}
\label{rem:diff}
The terminology \lq\lq{differential symmetry breaking operator}\rq\rq\
 in Definition \ref{defprop:regdiff}
 makes reasonable sense.  
In fact, 
 $T$ is a differential operator
 in the sense of Definition \ref{def:Peetre} 
 if and only if $\operatorname{Supp}T$
 is a singleton 
 in $H' \backslash G/H$
 (see \cite[Part I, Lemma 2.3]{KP1}).  
\end{rem}

The strategy of \cite{xtsbon}
 for the classification of {\it{all}} symmetry breaking operators
 for $(G,G')$ satisfying (PP)
 is to use the stratification 
 of $H'$-orbits in $G/H$ by the closure relation.  
To be more precise,
 the strategy is: 
\begin{itemize}
\item[$\bullet$]
to obtain all differential symmetry breaking operators, 
 which corresponds to the singleton
 in $H' \backslash G/H$, 
 or equivalently,
 to solve certain branching problems
 for generalized Verma modules
 (see Section \ref{subsec:finDiff} below)
 via the duality \eqref{eqn:Vdual}, 
\item[$\bullet$]
to construct and classify 
$
  \{T \in H(\lambda, \nu):
   \operatorname{Supp} T \subset \overline S\}
$
modulo $\{T \in H(\lambda, \nu):
   \operatorname{Supp} T \subset \partial S\}
$ for $S \in G' \backslash G/H$ inductively.   
\end{itemize}
The \lq\lq{F-method}\rq\rq\ \cite{xkhelgason, Eastwood60, KOSS, KP1} gives a conceptual 
 and a practical tool
 to construct differential symmetry breaking operators
in Step 1.  
The second step may involve analytic questions
 such as the possibility of an extension of 
 an $H'$-invariant distribution 
 on an $H'$-invariant subset of $G/H$
 satisfying a differential equation
 to an $H'$-invariant distribution solution on the whole of $G/H$
 (e.g., 
 \cite[Chapter 11, Sect. 4]{xtsbon}), 
 and an analytic continuation
 and residue calculus with respect to some natural parameter
 (e.g., 
 \cite[Chapters 8 and 12]{xtsbon}).

We expect 
 that the methods developed in \cite{xtsbon}
 for the classification
 of symmetry breaking operators
 for the pair $(G,G')=({\mathrm{O}}(n+1,1),{\mathrm{O}}(n,1))$
 would work for some other pairs $(G,G')$
such as those satisfying (PP)
 (see Theorem \ref{thm:PP}
 for the list), 
 or more strongly those satisfying (BB)
 (see Theorem \ref{thm:BB}
 for the list).

\subsection{Finiteness criterion
 for differential symmetry breaking operators}
\label{subsec:finDiff}

As we have seen in Theorem \ref{thm:fmXY}
 and Proposition \ref{prop:fmP}, 
 it is a considerably strong restriction
 on the $G'$-manifold $Y$
 for the space 
 $\invHom{G'}{C^{\infty}(X,{\mathcal{V}})}{C^{\infty}(Y,{\mathcal{W}})}$
 of symmetry breaking operators
 to be finite-dimensional, 
 which would be a substantial condition 
 for further study in Stages B and C
 of the branching problem.  
On the other hand, 
 if we consider only {\it{differential}} symmetry breaking operators,
 then it turns out that there are much broader settings
 for which the finite-multiplicity property
 (or even the multiplicity-free property) holds.  
The aim of this subsection
 is to formulate this property.

In order to be precise,
 we write $\invHom{G'}{C^{\infty}(X,{\mathcal{V}})}{C^{\infty}(Y,{\mathcal{W}})}$ for the space of continuous symmetry breaking operators, 
 and 
 $\operatorname{Diff}_{G'}(C^{\infty}(X,{\mathcal{V}}), C^{\infty}(Y,{\mathcal{W}}))$
 for that of differential symmetry breaking operators.  
Clearly we have
\begin{equation}
\label{eqn:DC}
\operatorname{Diff}_{G'}(C^{\infty}(X,{\mathcal{V}}), C^{\infty}(Y,{\mathcal{W}}))
\subset
\invHom{G'}{C^{\infty}(X,{\mathcal{V}})}{C^{\infty}(Y,{\mathcal{W}})}.  
\end{equation}

\vskip 1pc
We now consider the problem 
 analogous to Problem \ref{prob:tabc}
 by replacing the right-hand side of \eqref{eqn:DC}
 with the left-hand side.

For simplicity, 
we consider the case where ${\mathcal{V}} \to X$
 is a $G$-equivariant line bundle
 over a real flag manifold $G/P$, 
 and write ${\mathcal{L}}_{\lambda} \to X$
 for the line bundle
 associated to a one-dimensional representation $\lambda$
 of $P$.  
We use the same letter $\lambda$
 to denote the corresponding infinitesimal representation 
of the Lie algebra ${\mathfrak {p}}$, 
 and write $\lambda \gg 0$
 if $\langle \lambda|_{{\mathfrak{j}}}, \alpha \rangle \gg 0$
 for all $\alpha \in \Delta({\mathfrak {n}}^+, {\mathfrak {j}})$
 where ${\mathfrak {j}}$ is a Cartan subalgebra
 contained in the Levi part ${\mathfrak{l}}$
 of the parabolic subalgebra ${\mathfrak{p}}={\mathfrak{l}}+{\mathfrak{n}}^+$.

We say a parabolic subalgebra ${\mathfrak{p}}$
 of ${\mathfrak{g}}$
 is ${\mathfrak{g}}'$-{\it{compatible}} 
 if ${\mathfrak{p}}$ is defined
 as the sum of eigenspaces
 with nonnegative eigenvalues 
 for some hyperbolic element in ${\mathfrak{g}}'$.  
Then ${\mathfrak{p}}':= {\mathfrak{p}} \cap {\mathfrak{g}}'$
 is a parabolic subalgebra
 of ${\mathfrak{g}}'$
 and we have compatible Levi decompositions
 ${\mathfrak{p}} = {\mathfrak{l}} + {\mathfrak{n}}^+$
 and ${\mathfrak{p}}' = ({\mathfrak{l}} \cap {\mathfrak{g}}')
 + ({\mathfrak{n}}^+ \cap {\mathfrak{g}}')$.  
We are ready to state
 an answer to a question analogous 
 to Problem \ref{prob:tabc} (1) and (2)
 for {\it{differential}} symmetry breaking operators
 (cf. \cite{Eastwood60}).

\begin{theorem}
[local operators]
\label{thm:finVerma}
Let $G'$ be a reductive subgroup
 of a real reductive linear Lie group $G$, 
 $X=G/P$
 and $Y=G'/P'$
 where $P$ is a parabolic subgroup of $G$
 and $P'=P \cap G'$
 such that the parabolic subalgebra ${\mathfrak{p}}={\mathfrak{l}} + {\mathfrak{n}}^+$
 of ${\mathfrak{g}}$ is ${\mathfrak{g}}'$-compatible.   

\par\noindent
{\rm{(1)}}\enspace
{\rm{(finite multiplicity)}}
\enspace
For any finite-dimensional representations $V$
 and $W$ of the parabolic subgroups $P$ and $P'$, 
 respectively,
 we have
\[
\dim_{\mathbb{C}} \operatorname{Diff}_{G'}(C^{\infty}(X, {\mathcal{V}}), 
C^{\infty}(Y, {\mathcal{W}}))
 < \infty, 
\]
where ${\mathcal{V}} = G \times_P V$
 and ${\mathcal{W}} = G' \times_{P'}W$
 are equivariant vector bundles
 over $X$ and $Y$, 
 respectively.  
\par\noindent
{\rm{(2)}}\enspace
{\rm{(uniformly bounded multiplicity)}}\enspace
If $({\mathfrak{g}},{\mathfrak{g}}')$ is a symmetric pair
 and ${\mathfrak{n}}^+$ is abelian,
 then 
 for any finite-dimensional representation $V$ of $P$, 
\[
C_V:=\sup_W \dim_{\mathbb{C}} 
\operatorname{Diff}_{G'}(C^{\infty}(X,{\mathcal{V}}),
 C^{\infty}(Y,{\mathcal{W}}))
 <\infty.   
\]
Here $W$ runs over all finite-dimensional irreducible 
 representations of $P'$.  
Furthermore, 
 $C_V=1$ if $V$ is a one-dimensional representation 
$\lambda$ of $P$
 with $\lambda \gg 0$.  
\end{theorem}

\begin{qedproof}
The classical duality
 between Verma modules
 and principal series representations
 in the case $G=G'$
 (e.g., 
 \cite{xhaja})
 can be extended 
 to the context
 of the restriction
 of reductive groups $G\downarrow G'$, 
 and the following bijection holds
 (see \cite[Part I, Corollary 2.9]{KP1}):
\begin{multline}
\label{eqn:Vdual}
\operatorname{Hom}_{({\mathfrak{g}}',P')}
  (U({\mathfrak{g}}') \otimes_{U({\mathfrak{p}}')}W^{\vee},
   U({\mathfrak{g}}) \otimes_{U({\mathfrak{p}})}V^{\vee})
\\
\simeq
\operatorname{Diff}_{G'}(C^{\infty}(G/P, {\mathcal{V}}), 
C^{\infty}(G'/P', {\mathcal{W}})).  
\end{multline}
Here $(\lambda^{\vee}, V^{\vee})$
 denotes the contragredient representation 
 of $(\lambda,V)$.  
The right-hand side of \eqref{eqn:Vdual} concerns
 Case II (symmetry breaking) in Section \ref{sec:ep},  
 whereas the left-hand side of \eqref{eqn:Vdual} concerns
 Case I (embedding)
 in the BGG category ${\mathcal{O}}$.  
An analogous theory of discretely decomposable restriction
 in the Harish-Chandra category ${\mathcal{HC}}$
 (see Sections \ref{sec:ep} and \ref{sec:A})
 can be developed more easily and explicitly 
 in the BGG category ${\mathcal {O}}$, 
 which was done in \cite{K12}.  
In particular, 
 the ${\mathfrak {g}}'$-compatibility
 is a sufficient condition 
 for the \lq\lq{discrete decomposability}\rq\rq\
 of generalized Verma modules
 $U({\mathfrak {g}}) \otimes_{U({\mathfrak {p}})} F$
 when restricted to the reductive subalgebra ${\mathfrak {g}}'$.  
Thus the proof
 of Theorem \ref{thm:finVerma}
 is reduced to the next proposition.

\begin{proposition}
\label{prop:Verma}
Let ${\mathfrak {g}}'$ be a reductive subalgebra
 of ${\mathfrak {g}}$.  
Suppose that a parabolic subalgebra ${\mathfrak {p}}={\mathfrak {l}}+{\mathfrak {n}}^+$
 is ${\mathfrak {g}}'$-compatible.  
\par\noindent
{\rm{(1)}}\enspace
For any finite-dimensional ${\mathfrak{p}}$-module $F$ 
 and ${\mathfrak{p}}'$-module $F'$, 
\[
  \dim \operatorname{Hom}_{{\mathfrak{g}}'}
  (U({\mathfrak{g}}') \otimes_{U({\mathfrak{p}}')}F',
   U({\mathfrak{g}}) \otimes_{U({\mathfrak{p}})}F)
  <\infty.  
\] 
\par\noindent
{\rm{(2)}}\enspace
If $({\mathfrak {g}}, {\mathfrak {g}}')$ is a symmetric pair
 and ${\mathfrak {n}}^+$ is abelian,
 then 
\[
\sup_{F'} \dim 
\operatorname{Hom}_{{\mathfrak {g}}'}
(U({\mathfrak {g}}') \otimes_{U({\mathfrak {p}}')}F',
U({\mathfrak {g}}) \otimes_{U({\mathfrak {p}})} {\mathbb{C}}_{\lambda})
=1
\]
for any one-dimensional representation $\lambda$ 
 of ${\mathfrak {p}}$
 with $\lambda \ll 0$.  
Here the supremum is taken over all finite-dimensional
 simple ${\mathfrak {p}}'$-modules $F'$.  
\end{proposition}
\begin{qedproof}
(1)\enspace
The proof is parallel to \cite[Theorem 3.10]{K12}
which treated the case
 where $F$ and $F'$ are simple modules
 of $P$ and $P'$, 
respectively.  
\par\noindent
(2)\enspace
See \cite[Theorem 5.1]{K12}.  
\end{qedproof}

Hence Theorem \ref{thm:finVerma} is proved.  
\end{qedproof}

\begin{rem}
\label{rem:m2}
If we drop the assumption
 $\lambda \gg 0$
 in Theorem \ref{thm:finVerma} (2)
 or $\lambda \ll 0$ in Proposition \ref{prop:Verma} (2), 
 then the multiplicity-free statement may fail.  
In fact, 
 the computation in Section \ref{subsec:SL2}
 gives a counterexample
 where $({\mathfrak{g}}, {\mathfrak{g}}')
=
 ({\mathfrak{sl}}(2,{\mathbb{C}}) + {\mathfrak{sl}}(2,{\mathbb{C}}), 
  \operatorname{diag}({\mathfrak{sl}}(2,{\mathbb{C}})))
$;
 see Remark \ref{rem:2.6} (3).  
\end{rem}

\begin{rem}
\begin{enumerate}[{\rm (1)}]
\item 
{\rm{(Stage B)}}\enspace
In the setting of Proposition \ref{prop:Verma} (2), 
 Stage B in the branching problem 
 (finding explicit branching laws)
 have been studied 
 in \cite{K08, K12}
 in the BGG category ${\mathcal{O}}$
 generalizing earlier results 
 by Kostant and Schmid \cite{xschthe}.  
\item 
{\rm{(Stage C)}}\enspace
In the setting of Theorem \ref{thm:finVerma} (2), 
 one may wish to find 
 an explicit formula
 for the unique {\it{differential}} symmetry breaking operators.  
So far, 
 this has been done
 only in some special cases;
 see \cite{C75, vDijkPevzner}
 for the Rankin--Cohen bidifferential operator, 
 Juhl \cite{Juhl}
 in connection with conformal geometry, 
 and \cite{KOSS, KP1} 
 using the Fourier transform 
 (\lq\lq{F-method}\rq\rq\ in \cite{xkhelgason}).  
\end{enumerate}
\end{rem}

We end this subsection
 by applying Theorem \ref{thm:finVerma}
 and Theorem \ref{thm:A}
 to the reductive symmetric pair
 $(G,G')=({\mathrm{GL}}(n_1+n_2, {\mathbb{R}}),{\mathrm{GL}}(n_1, {\mathbb{R}})\times {\mathrm{GL}}(n_2, {\mathbb{R}}))$, 
 and observe a sharp contrast
 between differential and continuous symmetry breaking operators,
 i.e., 
 the left-hand and right-hand sides of \eqref{eqn:DC}, 
 respectively.  
\begin{examp}
\label{ex:GLpqr}
{\rm{
Let $n=n_1+n_2$ with $n_1$, $n_2 \ge 2$.  
Let $P$, $P'$ be minimal parabolic subgroups
 of 
\[
(G,G')=
({\mathrm{GL}}(n,{\mathbb{R}}), {\mathrm{GL}}(n_1,{\mathbb{R}}) \times {\mathrm{GL}}(n_2,{\mathbb{R}})), 
\]
respectively,
 and set $X=G/P$ and $Y=G'/P'$.  
Then:
\begin{enumerate}[{\rm (1)}]
\item 
For all finite-dimensional representations $V$ of $P$
 and $W$ of $P'$, 
\[
\dim_{\mathbb{C}} \operatorname{Diff}_{G'}
(\operatorname{Ind}_P^G(V)^{\infty}, \operatorname{Ind}_{P'}^{G'}(W)^{\infty})<\infty.  
\]
Furthermore
 if $V$ is a one-dimensional representation ${\mathbb{C}}_{\lambda}$
 with $\lambda \gg 0$
 in the notation 
of Theorem \ref{thm:finVerma}, 
then the above dimension is 0 or 1.  
\item 
For some finite-dimensional representations $V$ of $P$
 and $W$ of $P'$, 
\[
  \dim_{\mathbb{C}} \operatorname{Hom}_{G'}(\operatorname{Ind}_P^G(V)^{\infty}, \operatorname{Ind}_{P'}^{G'}(W)^{\infty})=\infty.  
\]
\end{enumerate}
}}
\end{examp}
\vskip 1pc

\subsection{Localness theorem in the holomorphic setting}
\label{subsec:local}
In the last example
 (Example \ref{ex:GLpqr}) 
  and also Theorem \ref{thm:Speh} in Section \ref{subsec:2.2}, 
 we have seen in the real setting
 that differential symmetry breaking operators
 are \lq\lq{very special}\rq\rq\
 among continuous symmetry breaking operators.  
In this subsection 
we explain the remarkable phenomenon
 in the holomorphic framework
 that any continuous symmetry breaking operator
 between two representations
 under certain special geometric settings 
 is given by a differential operator;
 see Observation \ref{obs:SL2} (1)
 for the ${\mathrm{SL}}(2,{\mathbb{R}})$ case.  
A general case is formulated
 in Theorem \ref{thm:local} below.  
The key idea of the proof is
 to use the theory
 of discretely decomposable restrictions
 \cite{xkInvent94, xkAnn98, xkInvent98}, 
 briefly explained in Section \ref{sec:A}.  
A conjectural statement
 is given in the next subsection.

Let $G \supset G'$ be real reductive linear Lie groups,
 $K \supset K'$ their maximal compact subgroups,
 and $G_{\mathbb{C}} \supset G_{\mathbb{C}}'$
 connected complex reductive Lie groups
 containing $G \supset G'$
 as real forms, 
respectively.  
The main assumption of this subsection is 
 that $X:=G/K$ and $Y:=G'/K'$ are Hermitian symmetric spaces. 
To be more precise,
 let $Q_{\mathbb{C}}$ and $Q_{\mathbb{C}}'$
 be parabolic subgroups
 of $G_{\mathbb{C}}$ and $G_{\mathbb{C}}'$
 with Levi subgroups
 $K_{\mathbb{C}}$ and $K_{\mathbb{C}}'$, 
 respectively,
 such that the following commutative diagram 
consists of holomorphic maps:
\begin{alignat}{4}
   Y=\,\, &\,\, G'/K' &&\,\, \subset  &&\,\, X =\,\, &&\,\, G/K
\notag
\\
\small{\text{Borel embedding}} &\quad\,\cap && && &&\quad\cap \,\,{\small{\text{Borel embedding}}}
\label{eqn:Borel}
\\
 & G_{\mathbb{C}}'/Q_{\mathbb{C}}' &&\,\, \subset  && && G_{\mathbb{C}}/Q_{\mathbb{C}}.  
\notag
\end{alignat}

\begin{theorem}
[{\cite[Part I]{KP1}}]
\label{thm:local}
Let ${\mathcal{V}} \to X$, ${\mathcal{W}} \to Y$
 be $G$-equivariant, $G'$-equivariant
 holomorphic vector bundles, 
respectively.  
\begin{enumerate}[{\rm (1)}]
\item 
{\rm{(localness theorem)}}
Any $G'$-homomorphism from 
${\mathcal{O}}(X, {\mathcal{V}})$ to ${\mathcal{O}}(Y, {\mathcal{W}})$
 is given by a holomorphic differential operator,
 in the sense of Definition \ref{def:Peetre}, 
with respect to a holomorphic embedding $Y \hookrightarrow X$.  
\end{enumerate}
We extend ${\mathcal {V}}$ and ${\mathcal {W}}$
 to holomorphic vector bundles
 over $G_{\mathbb{C}}/Q_{\mathbb{C}}$ and $G_{\mathbb{C}}'/Q_{\mathbb{C}}'$, 
respectively.  
\begin{enumerate}[{\rm (2)}]
\item 
{\rm{(extension theorem)}}
Any differential symmetry breaking operator in (1)
 defined on Hermitian symmetric spaces
 extends to
 a $G_{\mathbb{C}}'$-equivariant holomorphic differential operator
 ${\mathcal{O}}(G_{\mathbb{C}}/Q_{\mathbb{C}}, {\mathcal{V}})
 \to 
 {\mathcal{O}}(G_{\mathbb{C}}'/Q_{\mathbb{C}}', {\mathcal{W}})$
 with respect to a holomorphic map
 between the flag varieties
 $G_{\mathbb{C}}'/Q_{\mathbb{C}}' 
 \hookrightarrow G_{\mathbb{C}}/Q_{\mathbb{C}}$.  
\end{enumerate}
\end{theorem}
\begin{rem}
The representation $\pi$ on the Fr{\'e}chet space
 ${\mathcal{O}}(G/K, {\mathcal{V}})$ is a maximal globalization
 of the underlying $({\mathfrak {g}}, K)$-module
 $\pi_K$ in the sense of Schmid \cite{xschast}, 
 and contains some other globalizations
 having the same underlying $({\mathfrak {g}}, K)$-module
 $\pi_K$ 
 (e.g., 
 the Casselman--Wallach globalization $\pi^{\infty}$).  
One may ask whether an analogous statement holds
 if we replace $(\pi, {\mathcal{O}}(G/K, {\mathcal{V}}))$
 and $(\tau, {\mathcal{O}}(G'/K', {\mathcal{W}}))$
 by other globalizations
 such as $\pi^{\infty}$ and $\tau^{\infty}$.  
This question was raised by D. Vogan
 during the conference at MIT in May 2014.  
We gave an affirmative answer in \cite[Part I]{KP1}
 by proving 
 that the natural inclusions
\[
  \operatorname{Hom}_{G'}(\pi, \tau)
  \subset
 \operatorname{Hom}_{G'}(\pi^{\infty}, \tau^{\infty})
  \subset
  \operatorname{Hom}_{{\mathfrak {g}}', K'}(\pi_K, \tau_{K'})
\]
 are actually bijective in our setting.  
\end{rem}

\subsection{Localness conjecture for symmetry breaking operators
 on cohomologies}
\label{subsec:loc}

It might be natural 
 to ask a generalization of Theorem \ref{thm:local}
 to some other holomorphic settings, from 
holomorphic sections to Dolbeault cohomologies, 
and from highest weight modules
 to $A_{\mathfrak{q}}(\lambda)$ modules.  
 \begin{problem}
\label{prob:loc}
To what extent 
does the localness and extension theorem hold
 for symmetry breaking operators
 between Dolbeault cohomologies?
\end{problem}
In order to formulate the problem more precisely,
 we introduce the following assumption 
 on the pair $(G,G')$ 
 of real reductive groups:
\begin{equation}
\label{eqn:Kideal}
\text{$K$ has a normal subgroup
 of positive dimension
 which is contained in $K'$.}
\end{equation}
Here, 
 $K$ and $K'=K \cap G'$ are maximal compact subgroups
 of $G$ and $G'$, 
respectively, 
 as usual.  
We write $K^{(2)}$
 for the normal subgroup in \eqref{eqn:Kideal}, 
 ${\mathfrak {k}}_0^{(2)}$ for the corresponding Lie algebra, 
 and ${\mathfrak {k}}^{(2)}$ for its complexification.  
Then the assumption \eqref{eqn:Kideal} means 
 that we have direct sum decompositions
\[
   {\mathfrak {k}} = {\mathfrak {k}}^{(1)} \oplus {\mathfrak {k}}^{(2)}, 
    \qquad
   {\mathfrak {k}}' = {{\mathfrak {k}}'}^{(1)} \oplus {\mathfrak {k}}^{(2)} 
\]
for some ideals ${\mathfrak {k}}^{(1)}$ of ${\mathfrak {k}}$
 and ${{\mathfrak {k}}'}^{(1)}$ of ${\mathfrak {k}}'$, 
 respectively.  
The point here is 
 that ${\mathfrak {k}}^{(2)}$ is common to both ${\mathfrak {k}}$ and ${\mathfrak {k}}'$.

We take $H \in \sqrt{-1} {\mathfrak {k}}_0^{(2)}$, 
 define a $\theta$-stable parabolic subalgebra of ${\mathfrak {g}}$ by 
\[
   {\mathfrak {q}} \equiv {\mathfrak {q}}(H) = {\mathfrak {l}}+ {\mathfrak {u}}
\]
as the sum of eigenspaces
 of $\operatorname{ad}(H)$ with nonnegative eigenvalues, 
 and set $L:=G \cap Q_{\mathbb{C}}$
 where $Q_{\mathbb{C}}=N_{G_{\mathbb{C}}}({\mathfrak {q}})$
 is the parabolic subgroup of $G_{\mathbb{C}}$.  
Then $L$ is a reductive subgroup 
 of $G$
 with complexified Lie algebra ${\mathfrak{l}}$, 
 and we have an open embedding 
 $X:= G/L \subset G_{\mathbb{C}}/Q_{\mathbb{C}}$
 through which $G/L$ carries a complex structure.  
The same element $H$ defines 
 complex manifolds $Y:=G'/L' \subset G_{\mathbb{C}}'/Q_{\mathbb{C}}'$
 with the obvious notation.

In summary,
 we have the following geometry
 that generalizes \eqref{eqn:Borel}:
\begin{alignat*}{4}
   Y=\,\, &\,\, G'/L' &&\,\, \subset  &&\,\, X =\,\, &&\,\,\, G/L
\\
\small{\text{open}} &\,\quad\cap && && &&\quad\cap \,\,{\small{\text{open}}}
\\
 & G_{\mathbb{C}}'/Q_{\mathbb{C}}' &&\,\, \subset  && && G_{\mathbb{C}}/Q_{\mathbb{C}}.  
\end{alignat*}
It follows from the assumption \eqref{eqn:Kideal}
 that the compact manifold $K/L \cap K$ coincides with $K'/L' \cap K'$.  
Let $S$ denote the complex dimension 
of the complex compact manifolds
 $K/L \cap K \simeq K'/L' \cap K'$.

\begin{examp}
\label{ex:GLcpx}
{\rm{
\begin{enumerate}[{\rm (1)}]
\item 
{\rm{(Hermitian symmetric spaces)}}
Suppose that $K^{(2)}$ is abelian.  
Then $Y \subset X$ are Hermitian symmetric spaces,
 $S=0$, 
 and we obtain the geometric setting of Theorem \ref{thm:local}.  
\item 
$(G,G')=({\mathrm{U}}(p,q;{\mathbb{F}}), {\mathrm{U}}(p';{\mathbb{F}})\times {\mathrm{U}}(p'',q;{\mathbb{F}}))$
 with $p=p'+p''$
 for ${\mathbb{F}}={\mathbb{R}}$, ${\mathbb{C}}$, or ${\mathbb{H}}$, 
 and $K^{(2)}={\mathrm{U}}(q;{\mathbb{F}})$.  
Then neither $G/L$ nor $G'/L'$ is a Hermitian symmetric space
 but the assumption \eqref{eqn:Kideal} is satisfied.  
Thus the conjecture below applies.  
\end{enumerate}
}}
\end{examp}
For a finite-dimensional holomorphic representation $V$ of $Q_{\mathbb{C}}$, 
 we define a holomorphic vector bundle $G_{\mathbb{C}} \times_{Q_{\mathbb{C}}} V$
 over the generalized flag variety $G_{\mathbb{C}}/Q_{\mathbb{C}}$, 
and write ${\mathcal{V}}:=G \times_L V$
 for the $G$-equivariant holomorphic vector bundle over $X=G/L$
 as the restriction $(G_{\mathbb{C}} \times_{Q_{\mathbb{C}}} V)|_{G/L}$.  
Then the Dolbeault cohomology $H_{\bar \partial}^j(X,{\mathcal{V}})$
 naturally carries a Fr{\'e}chet topology
 by the closed range theorem
 of the $\bar \partial$-operator,
 and gives the maximal globalization
of the underlying $({\mathfrak {g}},K)$-modules,
 which are isomorphic to Zuckerman's derived functor modules
 ${\mathcal{R}}_{{\mathfrak {q}}}^j(V \otimes {\mathbb{C}}_{- \rho})$
\cite{Vogan81, xwong}.  
Similarly for $G'$, 
 given a finite-dimensional holomorphic representation $W$ of $Q_{\mathbb{C}}'$, 
 we form a $G'$-equivariant holomorphic vector bundle
 ${\mathcal{W}}:= G' \times_{L'} W$
 over $Y=G'/L'$
 and define a continuous representation of $G'$
 on the Dolbeault cohomologies
 $H_{\bar \partial}^j(Y,{\mathcal{W}})$.  
In this setting we have the discrete decomposability
 of the restriction
 by the general criterion
 (see Fact \ref{fact:5.5}).  
\begin{proposition}
\label{prop:uk}
The underlying $({\mathfrak {g}},K)$-modules $H_{\bar \partial}^j(X,{\mathcal{V}})_K$
 are $K'$-admissible.  
In particular,
 they are discretely decomposable as $({\mathfrak {g}}',K')$-modules.  
\end{proposition}

Explicit branching laws 
 in some special cases
 (in particular, 
 when $\dim V=1$)
of Example \ref{ex:GLcpx}
 (1) and (2)
 may be found in \cite{K08} and \cite{xgrwaII, xk:1}, 
respectively.  

We are now ready to formulate
 a possible extension
 of the localness and extension theorem 
 for holomorphic functions
 (Theorem \ref{thm:local})
 to Dolbeault cohomologies 
 that gives geometric realizations
 of Zuckerman's derived functor modules.  

\begin{conj}
\label{conj:local}
Suppose we are in the above setting, 
and let $V$ and $W$ be finite-dimensional representations of $Q_{\mathbb{C}}$ and $Q_{\mathbb{C}}'$, 
 respectively.  
\begin{enumerate}[{\rm (1)}]
\item 
{\rm{(localness theorem)}}\enspace
Any continuous $G'$-homomorphism 
\[
H_{\bar \partial}^S(X,{\mathcal{V}}) \to H_{\bar \partial}^S(Y,{\mathcal{W}})
\]
 is given by a holomorphic differential operator
 with respect to a holomorphic embedding $Y \hookrightarrow X$.  
\item 
{\rm{(extension theorem)}}
Any such operator in (1)
 defined on the open subsets
 $Y \subset X$
 of $G_{\mathbb{C}}'/Q_{\mathbb{C}}' \subset G_{\mathbb{C}}/Q_{\mathbb{C}}$, 
 respectively, 
 extends to a $G_{\mathbb{C}}'$-equivariant holomorhic differential operator
 with respect to a holomorphic map
 between the flag varieties $G_{\mathbb{C}}'/Q_{\mathbb{C}}' \hookrightarrow 
G_{\mathbb{C}}/ Q_{\mathbb{C}}$.  
\end{enumerate}
\end{conj}
The key ingredient of the proof
 of Theorem \ref{thm:local} for Hermitian symmetric spaces
 was the discrete decomposability
 of the restriction of the representation 
 (Fact \ref{fact:2.2} (2)).  
Proposition \ref{prop:uk} is a part
 of the evidence for Conjecture \ref{conj:local}
 in the general setting.

\runinhead{Acknowledgements}
The author thanks J.-L. Clerc, T. Kubo, T. Matsuki, B. \O rsted, T. Oshima, 
Y. Oshima, M. Pevzner, P. Somberg, V. Sou{\v c}ek, B. Speh
for their collaboration on the papers which are mentioned in this article.
This article is based on the lecture
 that the author delivered at the conference
 {\emph{Representations of reductive groups
in honor of David Vogan on his 60th birthday}}
at MIT, 19-23 May 2014. 
He would like to express his gratitude to the organizers,
Roman Bezrukavnikov, Pavel Etingof, George Lusztig, Monica Nevins, and 
Peter Trapa,
for their warm hospitality during the stimulating conference.

This work was partially supported by
        Grant-in-Aid for Scientific Research (A) (25247006), Japan
        Society for the Promotion of Science.

\end{document}